\documentclass[preprint,authoryear,12pt]{elsarticle}

\usepackage{amssymb}
\usepackage{amsmath,dsfont}
\usepackage{amsbsy}
\usepackage{amsthm}
\usepackage{amscd}
\usepackage{amsfonts}
\usepackage{subfigure}
\usepackage{psfrag}
\usepackage{nicefrac}
\usepackage{epsfig}
\usepackage{color}
\usepackage{algorithm}
\usepackage{algorithmic}

\usepackage{bm}
\usepackage{tikz}
\usetikzlibrary{patterns}

\usepackage{array,booktabs,tabularx,ragged2e,multirow} 

\graphicspath{    {./}    {./Bilder/}
}

\def\reff#1{\mbox{\rm(\ref{#1})}}
\newcommand{\vect}[1]{\boldsymbol{#1}}  
\newcommand{\tens}[1]{\vect{#1}}

\newcommand{\V}{\vect}
\newcommand{\T}{\vect}
\newcommand{\tr}{\operatorname{tr}}
\newcommand{\cof}{\operatorname{cof}}

\newcommand{\dev}{\operatorname{dev}}

\newcommand{\body}{\Omega}

\newcommand{\td}{\,\mathrm{d}}
\newcommand{\Teps}{\vect{\epsilon}}
\newcommand{\Tsigma}{\vect{\sigma}}
\newcommand{\Vchi}{\ensuremath{\mbox{\boldmath $\chi$}}}

\newcommand{\nablaq}[1]{\ensuremath{\nabla_{\hspace{-0.9mm}\mbox{\begin{scriptsize}$\mathrm{#1}$\end{scriptsize}}}}\hspace{0.9mm}}

\journal{Journal of the Mechanics and Physics of Solids}
\begin{document}

\begin{frontmatter}

\title 
{On the crack-driving force of phase-field models in linearized and finite elasticity
}

\author[siegen]{Carola Bilgen}
\author[siegen]{Kerstin Weinberg\corref{cor1}}

\cortext[cor1]{Corresponding author \\ E-mail address: kerstin.weinberg@uni-siegen.de (K. Weinberg).}

\address[siegen]{Lehrstuhl f\"ur Fest\"orpermechanik, Universit\"at Siegen, D-57076 Siegen, Germany}

\begin{abstract}
The phase-field approach to fracture has been proven to  be a mathematically sound and easy to implement method for computing crack propagation with arbitrary crack paths. 
Hereby   crack growth is driven by energy minimization resulting in a variational  crack-driving force. The definition of this force out of a  
tension-related  energy functional does, however, not always agree with the established failure criteria of fracture mechanics.
In this work different variational formulations for linear and finite elastic materials are discussed and   ad-hoc driving forces are presented which are motivated  by general fracture mechanical considerations. 
The superiority of the generalized approach is demonstrated by a series of numerical examples.


\end{abstract}

\end{frontmatter}

\section{Introduction} \label{Sec:Introduction}
A crack in a solid with domain $\body$ forms a new  surface of a priori unknown size, position and evolution.
In order to compute such  moving boundary problems, phase-field simulations of fracture have gained enormous popularity recently,
e.g., \cite{HenryLevine2004,Karma_etal2001,miehe2010phase,deborst2013,Borden_etal2012,miehe2014phase,sargado2018high}. 
%
The basic idea behind this approach is to introduce an additional field -the phase field- which indicates the material's state.  The continuous phase field  $z(\vect x,t)$ describes the solid state by $z=0$ and the broken state by $z=1$ and evolves  in space $\V x \in \body$ and time $t \in [0,\mathcal{T}]$.
Because the phase field is by definition a continuous field and the moving crack boundaries are 'smeared' over a small but finite zone of width $l_c$, the phase-field fracture approach constitutes a diffuse-interface formulation.

%
In computational mechanics the phase-field simulations of fracture have basically started with the seminal work of \cite{BourdinFrancfortMarigo2008}, 
who showed that the  diffuse-interface formulation converges in the limit $l_c\to 0$ to the sharp interface model of the classical Griffith theory of brittle fracture. Since then, 
numerous applications and enhancements have been examined. Among others \cite{teichtmeister2017phase,bleyer2018phase} expanded the phase-field model to anisotropic brittle fracture, \cite{ambati2015phase,ambati2016phase,kuhn2016phase} investigated ductile fracture in more detail and \cite{heider2017phase,ehlers2017phase} address to hydraulic fracture explicitly, to name some of them. 

%

Finite element simulations of fracture problems using the phase-field method have  mainly been motivated by the fact, that a rigorous connection between the phase-field and the sharp crack approach has been proven, regarding convergence of energy, energy release and evolution as the internal length $l_c$ vanishes, cf.~\cite{BourdinFrancfortMarigo2008,FreddyCarfagni2010,negri2013phase} and others. The link to Griffith's theory of fracture is usually made by means of global minimization arguments using  $\Gamma$-convergence theory for instance.  Such analyses obviously presume a variational setting, like it is employed in finite element methods to determine phase-field and deformation field by energy optimization. 

The propagation of a  crack   requires a state of tension at the crack tip. For mathematical analyses  this state is typically presumed  a priori, whereas for numerical simulations extra effort is needed to comply with this  physical constraint.  A common way to ensure a tensional state during crack growth is to split the body's energy into a tension-related  and a compressive energy functional    and to optimize for  the phase-field only the first. Thus, it results in a hybrid variational formulation, where the body's displacement field is obtained from the optimum of the full energy and the phase-field is obtained from the optimum of the tension-related part of the energy. Of course,  such a modification also allows for a mathematical rigorous treatment. The problem of the hybrid variational formulation is, however, the lacking connection to established theories of fracture mechanics.

The classical Griffith theory underlying the phase-field fracture approaches is a model which  is valid for ideal elastic and brittle materials under  tension. Its generalization, for example to ductile materials,
is not covered by the original theory,  	just as little as the common extensions of phase-field fracture models to materials showing  plastic  behavior.  Here 
we suggest to open the phase-field fracture approach to a wider range of applications by using physically motivated but not variationally determined crack-driving forces. Such driving forces do not minimize an energy potential because they are established ad-hoc, i.e. postulated from a  failure criterion for brittle (or ductile) materials.
In consequence, they allow to account for different types of fracture observed in the engineering practice. In this paper we show numerical examples computed
with different physically motivated crack-driving forces and discuss their validity. We  also show that some problems, like the compressive-split test of civil engineering, 
definitely require the use of a classical fracture mechanics strategy 
to obtain simulation results which match experimental observations.


The remaining of the paper is organized as follows:
At next  we will present the basic equations of the phase-field fracture approach    for linear and finite elasticity in Section~\ref{sec:basics}.
In Section~\ref{sec:drivingforce} the common energy splits for the linear and nonlinear   models are described and different definitions of their 
tensile components, which establish the variational crack-driving forces, are compared. Additionally to the  driving forces of the classical 
approach, a series of physically motivated ad-hoc driving forces are presented. Section~\ref{sec:numExamples} is devoted to the numerical examples. 
In the first part of this section we discuss various two-dimensional parametric studies, then follow two real-world   problems. 
We show  numerical simulations of a  Brazilian test and   a conchoidal fracture specimen
and discuss different crack-driving forces of linear and nonlinear elasticity.
A short summary closes the paper in Section~\ref{sec:conclusion}.

\section{The phase-field approach to fracture}\label{sec:basics}
Crack growth within a solid which of domain $\body  \subset\mathbb{R}^{3}$ and with boundary $\partial \body  \equiv \Gamma \subset\mathbb{R}^{2}$ corresponds to the creation of  new boundary surfaces $\Gamma(t)$. Therefore, the total potential energy of a homogenous but cracking solid is  composed of its material's energy  density $\Psi^\text{mat}$, and of a surface energy contribution from  growing crack boundaries.
\begin{align}\label{EpotBody}
  E  = \int_{\body}  \Psi^\text{mat}  \td \Omega + \int_{\Gamma(t)} {\mathcal{G}}_c \td \Gamma
\end{align}
The material's  energy  will be restricted within this investigation to an elastic free Helmholtz energy density, $\Psi^\text{mat}=\Psi^e$. The fracture-surface energy density $\mathcal{G}_c$ quantifies the material's resistance to cracking; for brittle fracture it corresponds to Griffith's critical energy release rate.
Both depend on the body's unknown displacement field $\vect{u}(\vect{x},t): \Omega\times [0,\mathcal{T}] \rightarrow \mathbb{R}^3$. In the course of solution $\V u(\V x,t)$ can be found when the energy attains a minimum.
Unfortunately, 
functional \reff{EpotBody} cannot be optimized in general and even an time-incremental technique is challenging because of the moving boundaries $\Gamma(t)$.

In the phase-field approach to fracture, the set of evolving crack boundaries is approximated by a surface-density function $\gamma=\gamma(z(\vect x,t))$ and
 \begin{align}\label{regularization}
  \int_{\Gamma(t)}   \ \text{d}\Gamma \approx \int_{\Omega}   \gamma(z(\vect{x},t)) \ \text{d}\Omega
 \end{align}
which allows to re-write the dissipative energy functional of a cracking solid and to formulate the optimization problem locally.
\begin{align}\label{Epotregularised0}
  E   =   \int_{\body} \left(  \Psi^e + {\mathcal{G}}_c \gamma   \right)\td \Omega \ \rightarrow \ \text{ optimum}
\end{align}
The elastic energy density is a function of   phase-field and  displacement, $\Psi^e(z,\V u)$.
By definition of the  surface-density function $\gamma(z)$, the fracture energy contribution $\mathcal{G}_c\gamma(z)$ is  only different from zero along cracks.

There are different ways to formulate the surface-density function, cf. \cite{Borden_etal2014,Hesch_etmany2015JCP}; 
here we will refer only to the second-order form
 \begin{align}\label{gamma2O}
 \gamma(z,\nabla z)=\frac{1}{2l_c}z^2+\frac{l_c}{2}|\nabla z|^2 \, .
 \end{align}
The first term results in a jump at the crack whereas the gradient term regularizes the discontinuity;  the length-scale parameter $l_c$ defines the width of the regularized, diffuse crack zone. Inserting the crack density function \eqref{gamma2O} in the total potential energy \reff{Epotregularised0} result in the well-known Ambrosio-Tortorelli functional of continuum damage mechanics, cf. \cite{AmbrosiTortorelli}. One main difference between the phase-field model and gradient damage models is that the order parameter can be just interpreted in the limits, $z=0$ indicates the intact solid material and $z=1$ the broken state. Values $z\in(0,1)$ have no physical meaning.

\subsection{Phase-field evolution equation}
The minimization of potential \reff{Epotregularised0} is in phase-field modeling typically restated as  a  gradient flow problem. This results in an evolution equation for the phase-field, $ \dot{z}(\vect{x},t) = \delta_z \Psi$, with $\Psi$ being the entire integrand of the potential \reff{Epotregularised0}. For balancing the units, a mobility parameter $M$ [$\mathrm{m^2/N\,s}$] is introduced which can also be understood as an inverse viscosity.
Collecting the driving forces and denoting their magnitude with $Y$  we may formulate the phase-field evolution as typical Allen-Cahn equation,
\begin{align}\label{eq:variationaldrivingforcedots}
    \dot{z}=MY(\vect{u},z)\,.
\end{align}
Normalizing the evolution equation gives a dimensionless effective driving force $\bar Y = \nicefrac{l_c}{\mathcal{G}_c}Y$,
\begin{align}\label{eq:variationaldrivingforcedotsDimless}
    \tau \dot{z}=\bar{Y}(\vect{u},z)\,,
\end{align}
where  $\tau = l_c/(M\mathcal{G}_c)$  is the retardation time [sec] of phase field evolution.

In the remaining of this paper we employ the normalized evolution equation \reff{eq:variationaldrivingforcedotsDimless}. The corresponding driving force can be written for a general, non-variational approach as
\begin{align}\label{eq:variationaldrivingforceYDimlessKLammer}
        \bar{Y} &=\bar{Y}^e +  l_c  \delta  \gamma \,   
\end{align}
where $\bar{Y}$ is postulated as the difference of an effective crack-driving force $\bar{Y}^e$ and the  geometrical crack resistance $l_c \delta\gamma $, i.e. the normalized regularizing term $\mathcal{G}_c \gamma$ in \reff{Epotregularised0}.
For the variational approach the generalized driving force can be formulated as
\begin{align}\label{eq:variationaldrivingforceDimlos}
    \bar{Y}=\delta \bar{\Psi} 
    = \delta \bar{\Psi}^e - (z + l_c^2 \triangle z)
\end{align}
with $\bar\Psi^e  = \nicefrac{l_c  }{ {\mathcal{G}}_c} \Psi^e(\V u,z) $.

In addition to that, in the cracked state the driving force has to vanish, $\bar Y(z=1) = 0$, such that the phase-field parameter stays limited. This is symbolized by the Macaulay brackets. Now, inserting the driving force $\bar Y$ into eq.~\reff{eq:variationaldrivingforcedotsDimless} the evolution equation can be re-written as
\begin{align}\label{eq:evolutionEqDimensionless}
\tau \dot z = \langle \bar Y^e + l_c\delta\gamma\rangle_+.
\end{align}

The viscous regularization is determined by the choice of retardation time $\tau$, or, equivalently, mobility parameter $M$ in \reff{eq:variationaldrivingforcedots}. If the retardation time $\tau= \nicefrac{l_c  }{ M\mathcal{G}_c}$ very small or the mobility is high, the phase-field changes suddenly and the effect of numerical regularization vanishes. If the mobility is too low and the retardation time too high, the phase-field evolution is suppressed and the structure stiffens artificially. Therefore, the retardation time has to be set is such a way, that the phase-field  is allowed to significantly evolve within one step of time discretization.  This is guaranteed, when it is in the order of magnitude of a time step, i.e.
\begin{align}\label{mobilityRule}
     \tau = O( \Delta t)
\end{align}
or $\tau \approx c \Delta t$ where a constant of $c=\nicefrac{1}{10} \dots 10$ still leads almost identical results. Translated into the mobility of eq.\reff{eq:variationaldrivingforcedots} this is $M \approx 2l_c/(c \Delta t\, \mathcal{G}_c)$.

\begin{figure}
 \centering
  \psfrag{force}{\small{$F_\text{max}$ [N]}}
  \psfrag{displacement}[cc][cc]{\small{$\bar u$ [mm]}}
 \includegraphics[width=0.45\linewidth]{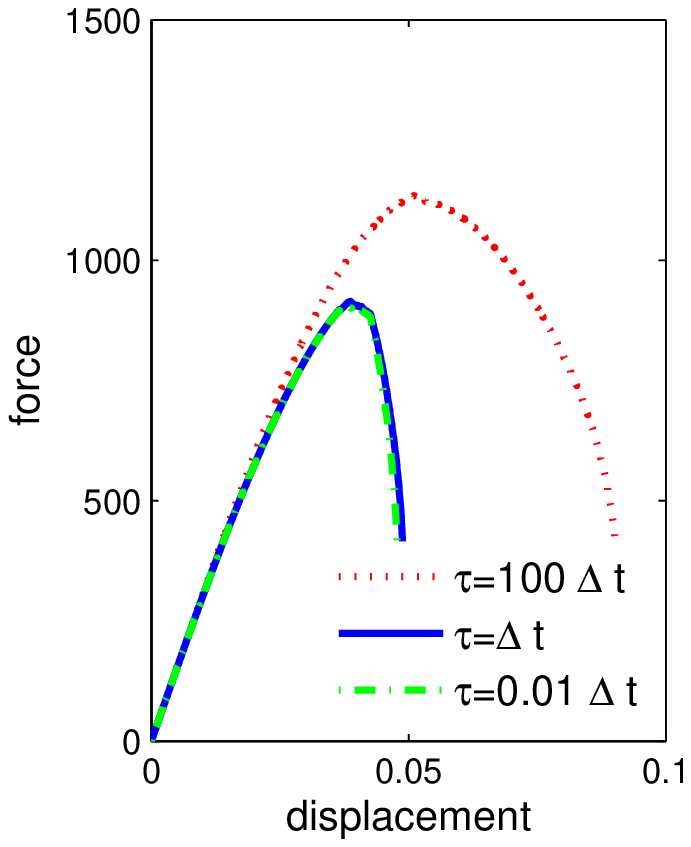}\qquad
 \includegraphics[width=0.45\textwidth]{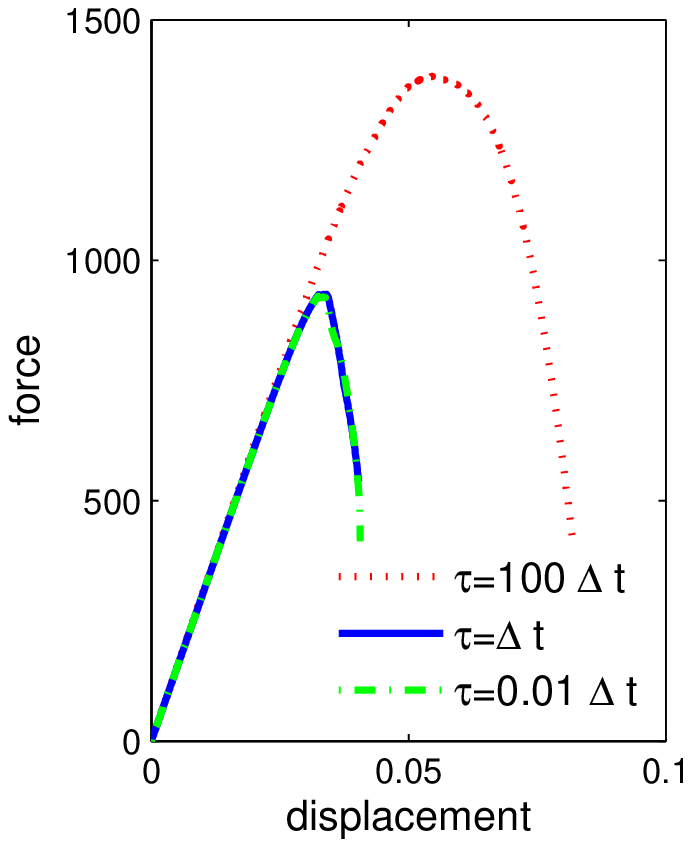}
 \caption{Load-displacement curves for a mode-I tension tension test computed with different relaxation times $\tau$ for an elastic energy driven crack growth (left) and for a principal stress driven crack growth (right)}
 \label{fig:mobilityGcSc}
\end{figure}

Figure~\ref{fig:mobilityGcSc} shows load-displacement curves for a mode-I tension test calculated with different values of $\tau$. The specific setup of the test is described in Sect.~\ref{sec:modeItensionTest}, the details of the crack propagation criteria are discussed below. Nonetheless we see, that for both crack-driving forces  a high retardation time results in a higher maximum loading before crack growth and, thus, an overrated load carrying capacity of the  structure.

Please note that the effect of viscous regularization holds  for the phase-field evolution in general. It is independent of the criterium of crack evolution, i.e. the specific crack-driving force on the right hand side of eq.~\reff{eq:variationaldrivingforcedotsDimless}. The strategy of \cite{miehe_etal_2015a,miehe2015phase}, where
an extra mobility parameter was introduced for the complementary energy as crack-driving force,  cannot be confirmed here.

\subsection{Linear elasticity}
In classical linear elasticity 
the strain-energy function is defined as a function of displacement $\V u$,
 \begin{align}\label{LinearStrainEnergyC}
  \Psi^e  = \frac{1}{2}\Teps(\vect{u}):\mathbb{C} :\Teps(\vect{u}),
 \end{align}
where $\Teps(\vect{u}) =\operatorname{sym} (\nabla \vect{u})$ is the strain tensor and  $\mathbb{C}$ is the Hookean material tensor which, for an isotropic material with Lam\'{e} constants $\lambda$ and $\mu$,  has the components   $C_{ijkl}=\lambda \delta _{ij}\delta _{kl}+\mu \left(\delta _{ik}\delta _{jl}+\delta_{il}\delta _{jk}\right)$.
Typically,  we formulate \reff{LinearStrainEnergyC} with the decomposition of the strain tensor into volumetric and deviatoric components,
$\Teps = \nicefrac13 \tr \Teps \T I + \dev \Teps$,
 \begin{align}\label{LinearStrainEnergy}
  \Psi^e  =  \frac12 \lambda \left( \tr \Teps \right)^2 + \mu \left| \Teps \right|^2 = \frac12 K\left( \tr \Teps  \right)^2 + \mu \left| \dev \Teps  \right|^2,
 \end{align}
where $K=\lambda + \nicefrac23 \mu$ is the  bulk modulus and $|\bullet|^2$ 
denotes the natural Frobenius norm. The elastic stresses are conjugate to the strains.
\begin{align}\label{stressesdWdeps}
    \Tsigma =\frac{\partial \Psi^e}{\partial \Teps } \stackrel{\phantom{\text{isotrop}}}{=}& \mathbb{C} :\epsilon(\vect{u})
    \\\nonumber
             \stackrel{\text{isotrop}}{=}& \lambda   \tr \Teps \T I  + 2\mu\Teps = K  \tr \Teps  \T I +2 \mu   \dev \Teps \,.
\end{align}

\subsection{Finite elasticity}
The general concept of Griffith's critical energy release rate and the corresponding potential energy \reff{Epotregularised0} is limited to brittle fracture but does not presume small deformations. The approach is also valid in the finite deformation regime with  deformation mapping $\Vchi(\V X,t):\Omega\times [0,\mathcal{T}] \rightarrow \mathbb{R}^3$ and $\V X+\V u(\V x,t)=\Vchi(\V X,t)$; the fields in capitals  refer to the initial configuration.
We define the deformation gradient and its determinant,
\begin{align}\label{defF}
\T{F}= \nablaq{\T{X}} {\Vchi} 
\quad\text{ and }\quad J = \det\tens{F}\,,
\end{align}
and also the right Cauchy-Green stress tensor $\tens{C}$ and its isochoric part,
\begin{align}\label{eq:CauchyGreen}
 \tens{C} = \tens{F}^{T}\tens{F}\quad\text{ and }\quad  \bar{\tens{C}}=J^{2/3}\bar{\tens{C}}\,.
\end{align}
To formulate the elastic strain energy functions of different material models it is useful to work with the  well-known principal invariants of the right Cauchy-Green deformation tensor for an isotropic elastic material, namely $I_1= \operatorname{tr}(\tens{C})$, $I_2 = \operatorname{tr}(\operatorname{cof} \tens{C})$ and $ I_3 = \operatorname{det}\tens{C}$.
They can be also be expressed with the principal eigenvalues $\lambda_a$, $a\in\{1,2,3\}$, e.g. $I_1=\lambda_1^2+\lambda_2^2+\lambda_3^2$.

Then, we define the elastic strain energy $ \Psi^e = \Psi^e(I_1, I_2, I_3)$, where an additional dependence on the phase-field enters in the presence of cracks. Exemplarily we state here the strain energy density function of a  classical Mooney-Rivlin material,
\begin{align}\label{MooneyRivlinAna}
 \Psi^e  = \frac{\mu}{2}[(I_1-3+k(I_2-3)]+\frac{K}{2}(J-1)^2
\end{align}
with initial shear modulus $\mu$,  second shear modulus $\mu_2$ and $k=\mu_2/\mu$; $K$ is again the bulk modulus.  Constitutive relation \reff{MooneyRivlinAna} implicitly requires an a priori incompressible material. For numerical computations (nearly) incompressible material behavior is commonly modeled by using a multiplicative decomposition of the deformation gradient into volume-changing and volume-preserving parts, $\tens{F} = J^{1/3}\bar{\tens{F}}$. For that reason, the adapted principal invariants are formulated with \reff{eq:CauchyGreen}$_2$ as
\begin{align}\label{eq:invariants}
 \bar I_1 &= \operatorname{tr}(\bar{\tens{C}}) = \det(\tens{F})^{-2/3}\tens{F}:\tens{F},\\
 \bar I_2 &= \operatorname{tr}(\text{cof} \bar{\tens{C}}) = \det(\tens{F})^{4/3}\tens{F}^{-T} : \tens{F}^{-T},\\
 \bar I_3 &= \text{det}\bar{\tens{C}} = J^2 = 1 \,,
\end{align}
and the elastic strain energy can be decomposed in its volumetric and isochoric contributions,
\begin{align}\label{PsiAdditiveSplitDevVol}
 \Psi^e = U(J) + \bar{\Psi}_0(\bar I_1,\bar I_2).
\end{align}
For an incompressible Mooney-Rivlin material the isochoric part of the strain energy density function \reff{PsiAdditiveSplitDevVol} reads
\begin{align}\label{MooneyRivlinNum}
    \bar{\Psi}_0(\bar{I}_1,\bar{I}_2)=\frac{\mu}{2}\left[\left(\bar{I}_1-3\right) +  k \left(\bar{I}_2^{3/2}-3\sqrt{3}\right)\right]\,.
\end{align}
This model deviates from the classical version \reff{MooneyRivlinAna} due to constraints following from the enforced polyconvexity. In \cite{HarNef03PGPH} it has been shown, that expressions $(\bar{I}_1-3)^{i} $ are convex in $\tens{F}$ for all $i\ge 1$, whereas convexity of the second invariant requires a form $(\bar{I}_2^{k}-3^k)^i $ with exponent $k\ge 3/2$.
The Neo-Hookean   model results from model \reff{MooneyRivlinNum} with $k= 0$. Concerning the volumetric contribution $U(J)$ there exist various approaches that have to fulfill a few conditions, cf. \cite{DollSchweizerhof2000}. In this paper the simplest version for the volumetric part is chosen as $U(J) = \frac{1}{2}K(J-1)^2$.

\section{Crack-driving forces}\label{sec:drivingforce}
The phase-field model introduced by \cite{FrancfortMarigo1998} refers to Griffith fracture, i.e. brittle material with crack growth proportional to the total elastic strain energy. However, for the computational simulation of fracture we need to account for the fact, that crack growth requires a state of tension whereas a compressive state --with same energy density-- does not result in crack growth.
This asymmetry of tension and compression prevents the usage of the same variational functional for both fields, the deformation $\V u$ and the phase-field $z$.

Instead, the variational formulation of the phase-field driving force \reff{eq:variationaldrivingforceDimlos} requires a split of the energy density into a tensile energy functional $\Psi^{e+}$ and a compressive energy  functional $\Psi^{e-}$ whereby only the first interacts with the  phase-field and drives the crack. To be a physically meaningful elastic split the identity
 \begin{align}\label{LinearSplit}
  \Psi^e(\V u,z) = \Psi^{e+}(\V u, z) + \Psi^{e-}(\V u)  \,.
 \end{align}
must hold, i.e. the  sum of tensile and compressive contribution must correspond to the total elastic energy. This is not always the case in phase-field fracture models as, in particular  for non-linear elasticity, the decomposition is somewhat arbitrarily.

\subsection{Variational crack-driving forces in linear elasticity}\label{sec:LinearSplit}
In linear elasticity the tension-compression asymmetry is typically modeled with  a decomposition of the  principal strains $\epsilon_a$, $a\in\{1,2,3\}$, into tensile and compressive components,
\begin{align}\label{MacaulayBrackets}
\left\langle \epsilon_a \right\rangle_+ = \frac12\left( \epsilon_a + |\epsilon_a|\right), \qquad
\left\langle \epsilon_a \right\rangle_- = \frac12\left( \epsilon_a - |\epsilon_a|\right), 
\end{align}
and a spectral representation of the   tensile and compressive strain tensors using the principal directions $\vect{n}_a$,
\begin{align}\label{Tepsplusminus}
    \Teps^{\pm}(\vect{u}) = \sum_{a=1}^{3} \left\langle \epsilon_a\right\rangle_\pm \vect{n}_a\otimes \vect{n}_a \,.
\end{align}
The corresponding linear-elastic strain energy densities can be stated as
\begin{align}\label{def:Psiplusminus}
     \Psi^{e\pm}_0 (\Teps^\pm) = \frac12 \Teps^\pm :\mathbb{C}:\Teps^\pm \,.
\end{align}
To see the validity of the decomposition \reff{def:Psiplusminus}, we consider the set of all possible elastic strain tensors $T$, i.e. $\Teps \in T$. The decomposed positive strain tensors \reff{Tepsplusminus} form a convex subset $\Teps^+ \in T^+$, $T^+\subset T$. The stresses calculated with the negative strain tensors,  $\mathbb{C}:\Teps^- \in T^-$, form a subset  which is dual to $T^+$. 
The latter follows from the general variational approach where we require $\Teps^+$ to minimize $\Teps-\Teps^+$ in the sense of the energy norm, cf. \cite{LiMarigo2017}. This results in an orthogonality condition, i.e. the products $\Teps^+ : \mathbb{C}\Teps^-$ vanish and it holds
\begin{align*}
    \Psi^{e}_0 = \frac{1}{2}(\Teps^{+}+\Teps^{-}):\mathbb{C}:(\Teps^{+}+\Teps^{-}) = \frac{1}{2}\Teps^{+}:\mathbb{C}:\Teps^{+} +
    \frac{1}{2}\Teps^{-}:\mathbb{C}:\Teps^{-} = \Psi^{e+}_0 + \Psi^{e-}_0 \,.
\end{align*}

Hence  the variation of $\Psi^{e+}_0$, which results from positive dilatation and positive principal strains only, can be used to drive the crack.
The phase-field enters the variational driving force by means of a degradation function $g(z):[0,1]\rightarrow[\varepsilon,1]$,
\begin{align}\label{PsigsPsiPlusMinus}
    \Psi^{e+} &= g(z) \Psi^{e+}_0
\end{align}
which  accounts for the loss of stiffness within the crack. Now the crack-driving force is
\begin{align}\label{crackdrivingforce}
    Y^e=\delta_z \Psi^{e+} \,.
\end{align}

The degradation function $g(z)$ is modeled with $g=1$ in the unbroken material and $g=0$ in the crack. The specific form of $g(z)$ is open to modeling as long as the conditions $g(0)=1$, $g(1)=0$, and reasonably $g'(1)=0$, are fulfilled, cf. \cite{WeinbergGAMMMitteilungen2016}. Here the common quadratic degradation function
\begin{align}\label{weightf}
 g(z) = (1-z)^2 +\varepsilon
\end{align}
is applied where the parameter $0<\varepsilon\ll 1$ is introduced to avoid numerical instabilities in the case of a fully broken state at $z=1$. For further discussions about the choice of degradation functions we refer to \cite{sargado2018high}.

From the variational consistent decomposition $\Psi^{e}_0 =\Psi^{e+}_0 + \Psi^{e-}_0$ we derive with \reff{stressesdWdeps}   the stresses
\begin{align*}
     \Tsigma &= g(z) \frac{\partial \Psi^{e+}_0}{\partial \Teps} + \frac{\partial \Psi^{e-}_0}{\partial \Teps}
\end{align*}
such that for the stresses the following decomposition remains
\begin{align}\label{SigmagsSigmaPlusMinus}
    \Tsigma &= g(z) \Tsigma^+_0 + \Tsigma^-_0 \,.
\end{align}
The stress tensor is by definition $\Tsigma=\partial \Psi/\partial \Teps$ where now, employing the energy function \reff{def:Psiplusminus}, derivatives of the decomposed strains $\Teps^{\pm}$ with respect to the total strain appear. These cannot be resolved in general and in order to apply decomposition \reff{def:Psiplusminus}, the unusual stress definition of $\Tsigma^\pm=\partial \Psi^\pm/\partial \Teps^\pm$ has to be used, cf~\cite{WeinbergDally2015}. Energy splits which allow to deduce the stress tensor from the classical definition go back to \cite{amor2009regularized,miehe2010phase} and make use of the common split into volumetric and deviatoric components to state
 \begin{align}\label{PsiKGPlus}
  \Psi^{e+}_0  =  \frac12 K \left\langle \tr \Teps \right\rangle_+^2 + \mu \dev \Teps^+ : \dev \Teps^+
 \end{align}
or
\begin{align}\label{PsilambdaPlus}
  \Psi^{e+}_0  =  \frac12 \lambda \left\langle \tr \Teps \right\rangle_+^2 + \mu \Teps^+ : \Teps^+ \,.
\end{align}
For $g(z)=1$ both decompositions result, for an arbitrary symmetric tensor $\Teps\in T$ and formulation \reff{stressesdWdeps}, in stresses $\Tsigma=\Tsigma^{+}+\Tsigma^{-}$. Thus, both decompositions constitute variational approaches.
 \begin{align}\label{SigmaKGPlus}
  \Tsigma^{\pm}  &=    K \left\langle \tr \Teps \right\rangle_\pm \tens{I} + 2 \mu \dev \Teps^\pm \,
 \\\label{SigmalambdaPlus}
  \Tsigma^{\pm}  &=    \lambda \left\langle \tr \Teps \right\rangle_\pm \tens{I} + 2 \mu \Teps^\pm \,
 \end{align}
Please note that the energy splits \reff{def:Psiplusminus}, \reff{PsiKGPlus} and \reff{PsilambdaPlus} are modeling assumptions to define a variational crack-driving force \reff{crackdrivingforce}. However, there are only two ways to base the crack-driving force  for $z$ on  the same variational principle like the displacement field $\V u$. The first way is to use the full strain-energy \reff{LinearStrainEnergyC}, which corresponds to the original model of \cite{FrancfortMarigo1998} but   generally does not match the physical restrictions. The second way would be to use a deviatoric model which
is common in gradient damage mechanics, where damage is induced by elastic strain tensors that have zero trace.
This, however, corresponds to a purely shear induced failure and  not  to brittle crack growth and so it also does not match the physics of the problem. Therefore, a hybrid variational approach, where $\V u$ follows from the optimun of $\Psi^{e}(\V u, z)$ and $z$ follows from the optimum of $\Psi^{e+}(\V u, z)$, is inevitable.

\begin{figure}
\centering
\psfrag{nurtr}{\small{$\frac12 K\left\langle \tr \epsilon\right\rangle^+ $}}
\psfrag{fulldev}{\small{$\frac12 K\left\langle \tr \epsilon\right\rangle_+^2 + \mu |\dev \epsilon |^2$}}
\psfrag{Kmu}{\small{$\frac12 K\left\langle \tr \epsilon\right\rangle_+^2 + \mu |\dev \epsilon |_+^2$}}
\psfrag{lambda}{\small{$\frac12 \lambda\left\langle \tr \epsilon\right\rangle_+^2 + \mu | \epsilon |_+^2$}}
\psfrag{eCe}{\small{$\frac12  \Teps^+ :\mathbb{C} :\Teps^+  $}}
\psfrag{eps}{$\epsilon_{11}$}
\psfrag{Gc/lc}{$\mathcal{G}_c/lc$}
\psfrag{psiplus}{$\Psi^+$}
     \includegraphics[width=0.78\linewidth]{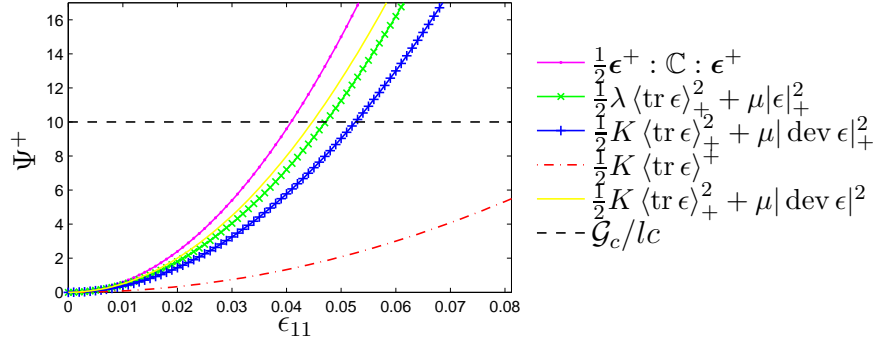}
	 \caption{Crack-driving energies $\Psi^{e+}$ for the   case of uniaxial tension in $1$-direction. The material parameters are: $E=10000$MPa, $\nu=0.25$, $\mathcal{G}_c=10$N/mm, $l_c=1\,$mm.}
	 \label{fig:1DZug_energies}
\end{figure}
\begin{figure}[b]
\centering
\psfrag{eps}{\tiny{$\epsilon_{11}$}}
\psfrag{psiplus}{\tiny{$\Psi^+$}}
    \includegraphics[width=0.32\linewidth]{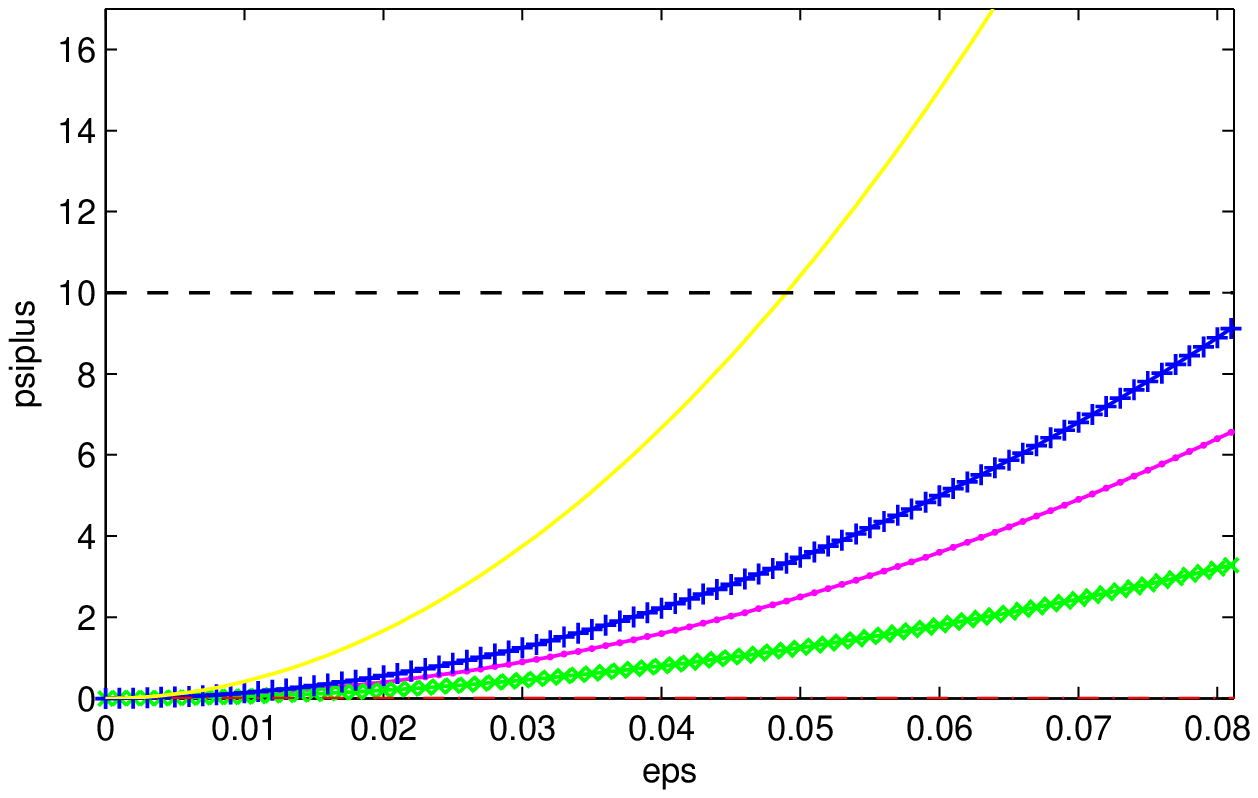}
    \includegraphics[width=0.32\linewidth]{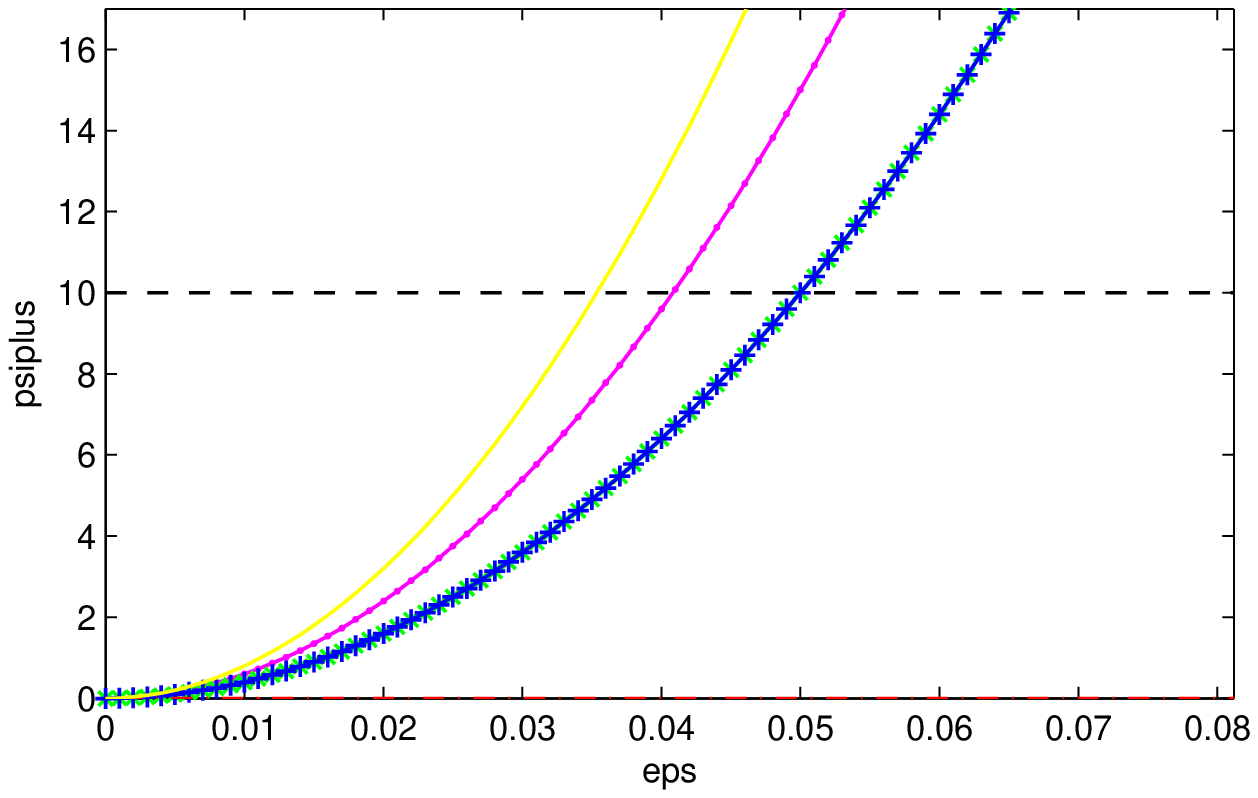}
    \includegraphics[width=0.32\linewidth]{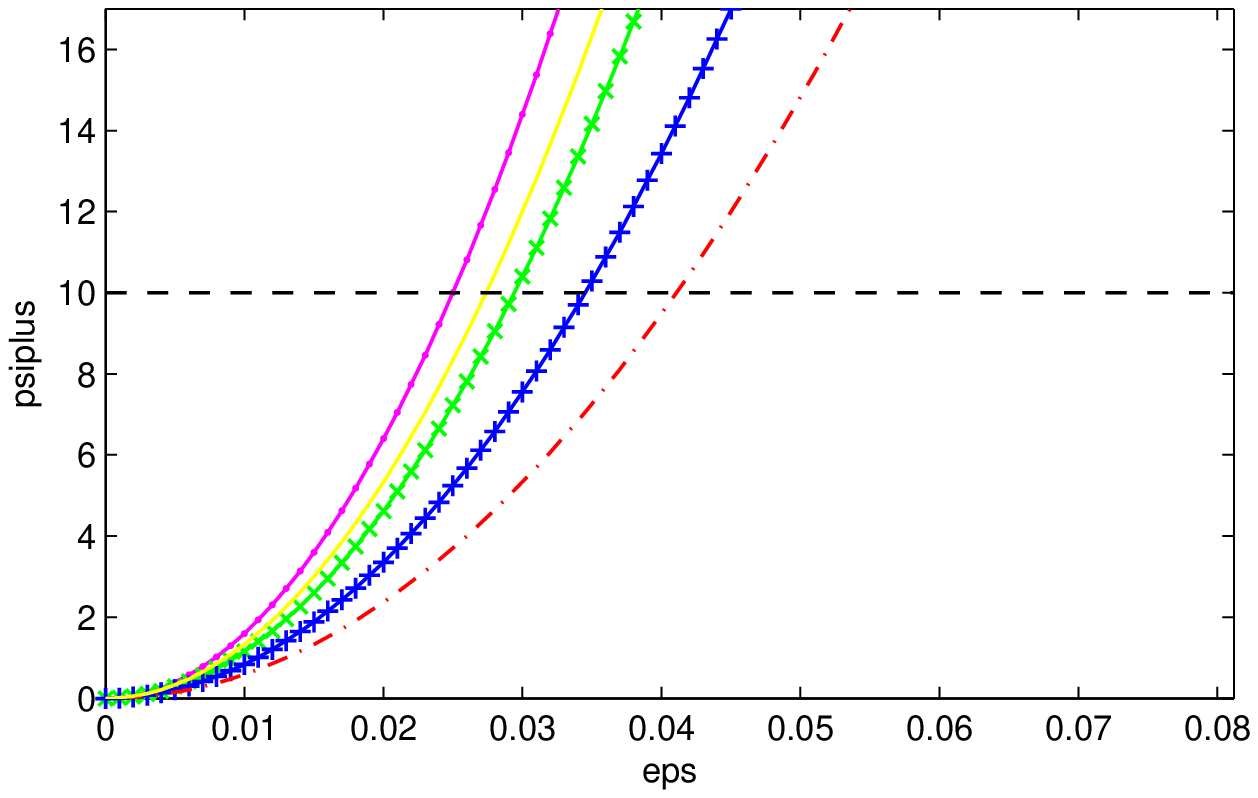}
	 \caption{Crack-driving energies $\Psi^{e+}$ for the   case of uniaxial compression (left), for pure shear (middle) and for      biaxial tension (right). The material parameters are: $E=10000\,$MPa, $\nu=0.25$, $\mathcal{G}_c=10\,$N/mm, $l_c=1\,$mm; the legend is the same as in Fig.~\ref{fig:2DZug_energies}}
	 \label{fig:2DZug_energies}
\end{figure}

%
The Figures~\ref{fig:1DZug_energies} and \ref{fig:2DZug_energies} illustrate the variational energy split for simple uniaxial loading cases. It can be seen,  that different decompositions result in different values for the crack-driving force.  Only for pure shear the results of \reff{PsiKGPlus} and \reff{PsilambdaPlus} lay on top of each other. Nonetheless a quantitative comparison is difficult because there is no 'reference energy split' in  fracture mechanics. Moreover, the maximum load depends on the choice of $\mathcal{G}_c$,  which is a material parameter but also on $l_c$ which is basically prescribed by numerical conditions.

\subsection{Variational crack-driving forces in finite elasticity}
Also for finite deformations the tension and compression asymmetry has to be regarded during the numerical simulations. Here we outline two different approaches. 
Both models will be introduced just shortly and for further information we refer to previous works, \cite{WeinbergHesch2015,Hesch_etmany2017,Thomas_etal_2018_GAMM}.

First we use the analogy to the small strains theory and formulate a spectral decomposition of the deformation gradient,
\begin{align}
   \T F  = \sum_{a=1}^{3} \lambda_a \vect{n}_a\otimes \vect{N}_a \,,
\end{align}
where $\lambda_a$, $a\in\{1,2,3\}$, denote the principal stretches;  $\vect{N}_a$ and $\vect{n}_a$ are the corresponding principal directions referring to the initial and current configuration, respectively.   The principal stretches $\lambda_a$ are   decomposed into tensile and compressive components,
\begin{align}\label{lambdaplusoderminus}
    \lambda_a^\pm =1+ \frac12\left( (\lambda_a-1)\pm |\lambda_a-1| \right) = 1+\langle \lambda_a-1 \rangle_\pm
\end{align}
and it holds  
\begin{align}\label{lambdaplusminus}
    \lambda_a =\lambda_a^+ \lambda_a^-\,,
\end{align}
%
%
in such a way, that only the elastic stretches $\lambda_a^e$ are driving the crack growth. They are weighed with the phase field to separate into elastic (tensile) and inelastic (compressive/cracked) components,
\begin{align}\label{lambdaei}
    \lambda_a^e =\left(\lambda_a^{+}\right)^{1-z} , \qquad
    \lambda_a^i =\left(\lambda_a^{+}\right)^{z}\lambda_a^-.
\end{align}
The exponent $1-z$ controls the tensile weakening during cracking, i.e., $\lambda_a^e = 1$ 
holds for a fully damaged material point with  $z=1$. 
Because of the multiplicative form \reff{lambdaei} the corresponding elastic strain energy density does not have an additive structure
in the sense of identity~\reff{LinearSplit},
\begin{align}
    \Psi^{e+} = \Psi^e(\lambda_1^e,\lambda_2^e,\lambda_3^e,z),
\end{align}
but the stresses can be derived in the usual way.  Specifically, the first Piola-Kirchhoff stress tensor is given by
\begin{align}\label{eq:ableitungPauslambdae}
    \tens{P}= \sum_{a=1}^{3} \frac{\partial\Psi^e}{\partial \lambda_a^e}
    \frac{\partial \lambda_a^e}{\partial\lambda_a }  \vect{n}_a\otimes \vect{N}_a \,.
\end{align}
A problem of this decomposition is its partial lack of convexity (as a results of exponent $0<z-1<1$) which might lead to 
problems during the numerical simulation, cf. \cite{Hesch_etmany2017}. Modified degradation functions can reduce the problem but do not solve it.

Thus, we introduce at next an alternative decomposition which is numerically stable, cf. \cite{Hesch_etmany2017}. This approach starts directly with the decomposition of the  invariants \reff{eq:invariants}. For a general   material model of finite elasticity we split the compressible component $J$ and the incompressible invariants $\bar{I}_1, \bar{I}_2$ into tensile and compressive parts to formulate
\begin{align}\label{SplitInv}
  \bar{I}_{1}^{\pm} &= 3+  J^{-2/3}  \langle\tens{F} : \tens{F}  -3 \rangle_\pm  \,\\
  \bar{I}_{2}^{\pm} &= 3+  J^{-4/3}  \langle\operatorname{cof}\tens{F} : \operatorname{cof}\tens{F}  -3 \rangle_\pm  \, \\
  J^{\pm} &=  1+ \langle\det \vect{F} -1 \rangle_\pm \,.
 \end{align}
The strain energy density \reff{PsiAdditiveSplitDevVol} has now the form
\begin{align}\label{SplitInvPsiplusminus}
    \Psi^e&= g_0(z)U(J^+) + U(J^-) + g_1(z)\Psi_0\left(\bar{I}_{1}^+,\bar{I}_{2}^+\right) + \Psi_0\left(\bar{I}_{1}^-,\bar{I}_{2}^-\right) \\\nonumber
    &\equiv g_0(z)U^+_0 + U^-_0 + g_1(z)\Psi_0^{+}  + \Psi^{-}_0\,,
\end{align}
where the degradation functions $g_0$, $g_1$ can differ from each other, in general. Our numerical experience has shown, however, that different functions are arbitrary and so we set $g_0=g_1$, which also enables a sound mathematical analysis of the positive part of \reff{SplitInvPsiplusminus}, cf. \cite{Thomas_etal_2018_GAMM}. Energy function~\reff{SplitInvPsiplusminus} can also be  formulated with arguments $J$, $\tr(\T F)$ and $\cof(\T F)$
and, thus, it is polyconvex in the mechanical deformation for any value of  phase-field $z$.
The first Piola-Kirchhoff stress tensor results from the derivative of the strain energy function with respect to the deformation gradient $\tens{F}$ as follows
\begin{align}
    \tens{P}  =\frac{\partial\Psi^{e}}{\partial\tens{F}} =  \frac{\partial g_0U^{+}}{\partial\tens{F}}
    + \frac{\partial U^-}{\partial\tens{F}} + \frac{\partial g_1\Psi^+_0}{\partial\tens{F}}
    + \frac{\partial \Psi^-_0}{\partial\tens{F}}.
\end{align}

Again, the crack-driving force is given by the variational derivative of the strain energy function \reff{SplitInvPsiplusminus}, i.e. $Y^e = \delta_z \Psi^e$.

\subsection{Fracture mechanically motivated crack-driving forces}
By now we have formulated the  crack-driving forces in a variational way, and the variational principle has been modified to account for the tension-compression asymmetry. However, there is no general method to  split the energy into tensile and compressive components. All decompositions are somewhat arbitrarily and at best justified  by mathematical considerations.
In particular, it is not possible to account for different fracture mechanical failure types, like rupture  and tear, through a tension-compression energy split. Therefore we chose here an alternative approach and model ad-hoc the crack-driving force $\bar{Y}^e$ in eq.~\reff{eq:variationaldrivingforceYDimlessKLammer}, motivated by specific fracture mechanical models.

To this end we first define the principal stresses $\sigma_a$, $a\in\{1,2,3\}$.  They are determined by solving the eigenvalue problem of Cauchy stress tensor $\Tsigma$ and will be ordered here, with $\sigma_I=\max(\sigma_a)$, $\sigma_{III}=\min(\sigma_a)$, and
\begin{align}\label{sortierteHauptspannungen}
    \sigma_{I} \geq \sigma_{II} \geq \sigma_{III}
\end{align}
so that with $\sigma_{I}(\V x,t) $ the maximum normal stress field is known. For the principal strains we refer to Sect.~\ref{sec:LinearSplit}.

\subsubsection{Maximum principal stress}

In classical mechanics failure is expected when the maximum tensile stress exceeds a critical value. Such considerations go back to Galilei, Navier and Rankine
and others, cf. \cite{GrossSelig2011}, and suggest to model crack growth to be driven by the maximum normal stress. As a threshold, given by  a material parameter of critical strength,  we introduce the decohesion stress $\sigma_c$. For brittle materials this decohesion stress is usually identified with the material's tensile resistance $R_m^t$, reduced values can be justified for small-scale plasticity in the vicinity of the crack tip.

\begin{figure}[b]\centering
 \psfrag{A}[cc]{a)}
 \psfrag{B}[cc]{b)}
 \psfrag{C}[cc]{c)}
 \psfrag{S1}[cc]{$\sigma_{1}$}
 \psfrag{S2}[cc]{$\sigma_{2}$}
 \psfrag{S3}[cc]{$\sigma_{3}$}
  \includegraphics[width=0.99\linewidth]{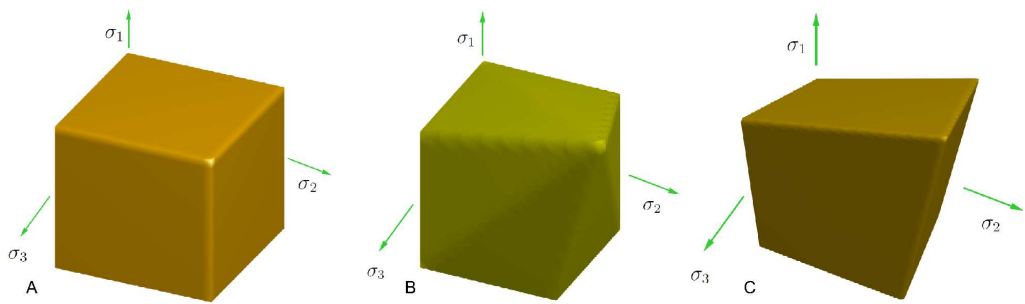}
   \caption{Elastic limit stresses for (a) the Rankine model , (b) the Mohr-Coulomb model and (c) the Beltrami model of the crack-driving force}
   \label{fig:bruchflachenI_III-IV}
\end{figure}

Now, assuming  that the crack starts   propagating  when the normal stress exceeds the decohesion stress $\sigma_c$, we define the Rankine model for the
crack-driving force
\begin{align}\label{Y0ausHauptspannung}
    \bar{Y}^e =\left\langle\frac{\sigma_I}{\sigma_c}-1\right\rangle_+.
\end{align}
This model corresponds to the classical failure criterion of  fracture mechanics and is displayed in the principal stress space
in Fig.~\ref{fig:bruchflachenI_III-IV}.

\subsubsection{Maximum shear stress}
While the normal stress is responsible for the crack growth in mode I direction, it might be possible that high shear stresses influence the crack propagation in mode II and III direction. Therefore, we state the  
maximum shear stress
\begin{align}\label{MaxSchubspannung}
    \tau_I= \frac{1}{2}\max\left(|\sigma_1-\sigma_2|,|\sigma_1-\sigma_3|,|\sigma_2-\sigma_3|\right)
    \stackrel{\text{eq.}(\ref{sortierteHauptspannungen})}{=}\frac12 \left(\sigma_{I}-\sigma_{III}\right)\,,
\end{align}
and formulate the crack-driving force in the sense of a Tresca model.
\begin{align}\label{Y0ausSchubspannung}
    \bar{Y}^e = \left\langle\frac{\sigma_I-\sigma_{III}}{\sigma_c}-1\right\rangle_+
\end{align}
The crack-driving force may also be formulated with a critical shear stress $\tau_c$,
\begin{align}\label{Y0ausSchubspannungTau}
    \bar{Y}^e = \left\langle\frac{\tau_I}{\tau_c}-1\right\rangle_+ \,
\end{align}
and we may avoid the eigenvalue decomposition by calculating the maximum shear stress from the deviatoric stress components.
\begin{align}\label{SchubspannungDeviator}
    \tau_I  
    = \sqrt{\frac38 \dev \Tsigma : \dev \Tsigma }
\end{align}
It remains to remark that, in general, the maximum shear stresses are responsible for plastic deformations and are not experimentally confirmed   as crack growth criterion.

\subsubsection{Extended principal stress}
Because a multiaxial loading may affect crack propagation, an extended principal stress criterion is introduced which also takes the hydrostatic stress state into account. First we define the mean stress
\begin{align}
    \sigma_m=\frac{1}{3}\tr(\Tsigma) = \frac{1}{3}\left(\sigma_{I}+\sigma_{II}+\sigma_{III} \right)
\end{align}
and formulate the crack-driving force with the indicator function $\mathds{1}$
\begin{align}\label{Y0erweiterteNormalspannungen}
    \bar{Y}^e = \mathds{1}_{\{\sigma_m>0\}}\left\langle\frac{\sigma_I}{\sigma_c}-1\right\rangle_+
\end{align}
In this case the crack growth is driven by the maximum principal stress but only if the material is under hydrostatic tension.

\subsubsection{Frictional shear stress}
As an example for a more involved failure criterion we chose here the Mohr-Coloumb frictional shear stress. It is developed for materials where the tensile material resistance $R_m^t$  differs from the compressive resistance $R_m^c$. Typically we have $R_m^t<R_m^c$ and
an interplay of tension and shear induces, e.g.  cracks in concrete  and shear zones in rocks.
Making use of the sorted principal stresses \reff{sortierteHauptspannungen} and setting $\sigma_c=R_m^t$ the Mohr-Coloumb crack-driving force is:
\begin{align}\label{Y0nachMohrColoumb}
     \bar{Y}^e &= \left\langle \frac{\sigma_I}{R_m^t} -   \frac{\sigma_{III}}{R_m^c} -1\right\rangle_+
\end{align}
Alternatively, the effective Mohr-Coloumb stress can be written with the relation of tensile and compressive strength, $m={R_m^c}/{R_m^t}$, as
\begin{align}\label{stressnachMohrColoumb}
   \sigma_\text{eff} & = (m+1)\tau_I + (m-1)\sigma_m= m \sigma_I -  \sigma_{III}
\end{align}
and then the driving force \reff{Y0nachMohrColoumb} reads simply
\begin{align}\label{Y0nachMohrColoumb2}
    \bar{Y}^e &= \left\langle \frac{\sigma_\text{eff}}{\sigma_c} -   1\right\rangle_+\,.
\end{align}
Please note that the Mohr-Coloumb crack-driving force \reff{Y0nachMohrColoumb2} is here understood as a place-holder for any alternative crack growth model. Other formulations of the effective stress are possible.  There are numerous failure hypotheses, like the Drucker-Prager version, the Cam-Clay model or anisotropic criteria, which might be employed as well.

\subsubsection{Maximum principal strain}
By now we regard the stresses to induce cracks. However, also the strain state can be used to define a failure criterion, as well as a combination of both.
According to the classical hypotheses of St. Vernant and Beltrami, failure occurs when the maximum principle strain reaches a material specific critical value.
Thus we define the crack-driving force using the maximum principal strain as
\begin{align}\label{Y0ausHauptdehnung}
     \bar{Y} = \left\langle \frac{\epsilon_I}{\epsilon_c}-1\right\rangle_+
\end{align}
with the maximum principal strain $\epsilon_I=\max(\epsilon_a)$, $a\in\{1,2,3\}$,  and the critical decohesion strain $\epsilon_c$.

Let us remark, that ad-hoc crack-driving force models can be used in both, linear and finite, elasticity. Whereas the models (\ref{Y0ausHauptspannung}-\ref{Y0nachMohrColoumb2}) apply directly, model \reff{Y0ausHauptdehnung} has to be reformulated in principle stretches for finite deformations,
\begin{align}\label{Y0ausHauptdehnunglambda}
     \bar{Y} = \left\langle \frac{\lambda_I}{\lambda_c}-1\right\rangle_+\,.
\end{align}

\begin{figure}
\centering
    \includegraphics[width=0.75\linewidth]{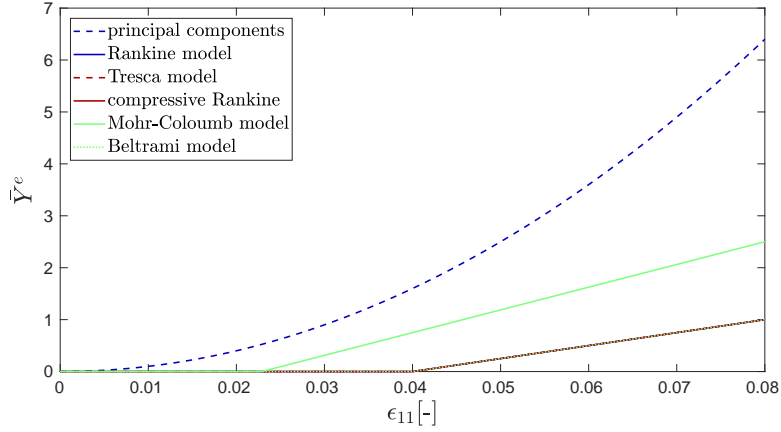}
	 \caption{Various crack-driving forces $Y^e$ in uniaxial strain. The material parameters are: $E=10000\,$MPa, $\nu=0.25$, $\mathcal{G}_c=10\,$N/mm, $l_c=1\,$mm.}
	 \label{fig:risstriebkraft}
\end{figure}%

\bigskip
In Figure~\ref{fig:risstriebkraft} all driving forces are plotted for the  case of simple uniaxial tension. It can be seen that the driving force in the energy based criterion takes immediate effect whereas the   stress-strain based  crack-driving forces start from a threshold. Also, the ad-hoc crack-driving forces are smaller than the positive energy. This may require an adapted relaxation time $\tau$ in \reff{eq:variationaldrivingforcedotsDimless} or, correspondingly, a somewhat higher mobility parameter. The driving forces of maximum principal stress, maximum shear stress, extended principal stress and maximum principal strain   coincide in the displayed loading case but not in general. Moreover, the slope of the driving force of  the Mohr-Coloumb model depends on the choice of material parameter $m$.

In Table~\ref{tab:CrackDrivingForces} all stated driving forces are summarized. Here we also define the short name of the models which will be referred  to in our sample computations.

\begin{small}
\begin{table}[htb]\label{tab:CrackDrivingForces}
\caption{Crack-driving forces for the phase-field evolution $\tau \dot{z} = \bar{Y}$ with $\bar{Y}=\bar{Y}^e+l_c \gamma$.} 
\label{tab:CrackDrivingForces} 
\renewcommand{\arraystretch}{1.2}%
\begin{tabular}{@{}llr}
\hline
  model & $\bar{Y}$ & short name\\ \hline
  elastic energy, $\delta_z \Psi^e$ & $\bar{Y}^e=   g'\left(\frac{\lambda}2  (\tr \Teps)^2 + \mu \Teps : \Teps \right)$ & Griffith model \\
  energy split \reff{def:Psiplusminus}, $\delta_z \Psi^e$ & $\bar{Y}^e=   g'\left( \Tsigma^+:  \Teps^+ \right)$ & principal components \\
   energy split \reff{PsilambdaPlus}, $\delta_z \Psi^+$ & $\bar{Y}^e=   g'\left(\frac{\lambda}2  \langle \tr \Teps\rangle_+^2 + \mu (\Teps^+)^2 \right)$ &$\lambda$-$\mu$ energy split \\
  energy split \reff{PsiKGPlus}, $\delta_z \Psi^+$ & $\bar{Y}^e=   g'\left(\frac{K}2 \langle \tr \Teps\rangle_+^2 + \mu \dev (\Teps^+)^2\right)$ & $K$-$\mu$ energy split \\
  maximum shear stress    & $\bar{Y}^e= \left\langle\frac{\tau_I}{\tau_c}-1\right\rangle_+$& Tresca model \\
  maximum principal stress    & $\bar{Y}^e=\left\langle\frac{\sigma_I}{\sigma_c}-1\right\rangle_+$ & Rankine model \\
  extended principal stresses    & $\bar{Y}^e= \mathds{1}_{\{\sigma_m>0\}}\left\langle\frac{\sigma_I}{\sigma_c}-1\right\rangle_+$ & compressive Rankine  \\
%
%
  frictional shear stress    & $\bar{Y}^e=  \left \langle \frac{\sigma_I}{R_m^t}-\frac{\sigma_{III}}{R_m^p}-1\right\rangle_+$ & Mohr-Coulomb model \\
  maximum principal strain    & $\bar{Y}^e= \left\langle\frac{\epsilon_I}{\epsilon_c}-1\right\rangle_+$  & Beltrami model \\
 \hline
\end{tabular}%
\\[2pt]
\end{table}%
\end{small}

\section{Numerical examples}\label{sec:numExamples}
\subsection{Plane mode-I tension and mode-II shear tests}\label{sec:modeItensionTest}
\begin{figure} 
 \centering
  \psfrag{U}[cc]{$\bar{u}$}
 \includegraphics[width= 0.3\linewidth]{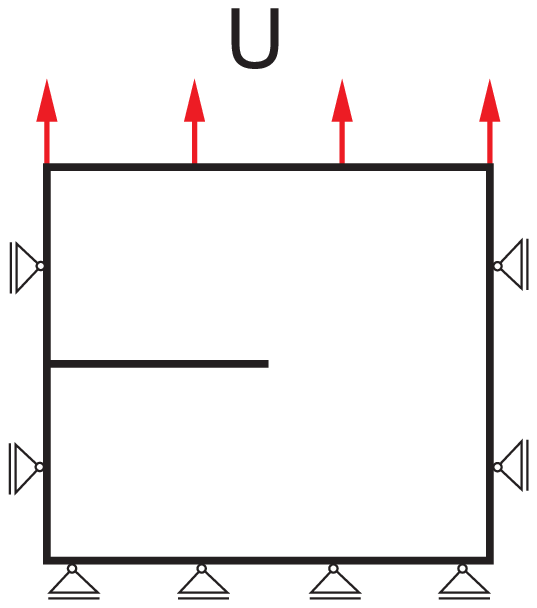}
 \qquad \includegraphics[width=0.3\textwidth]{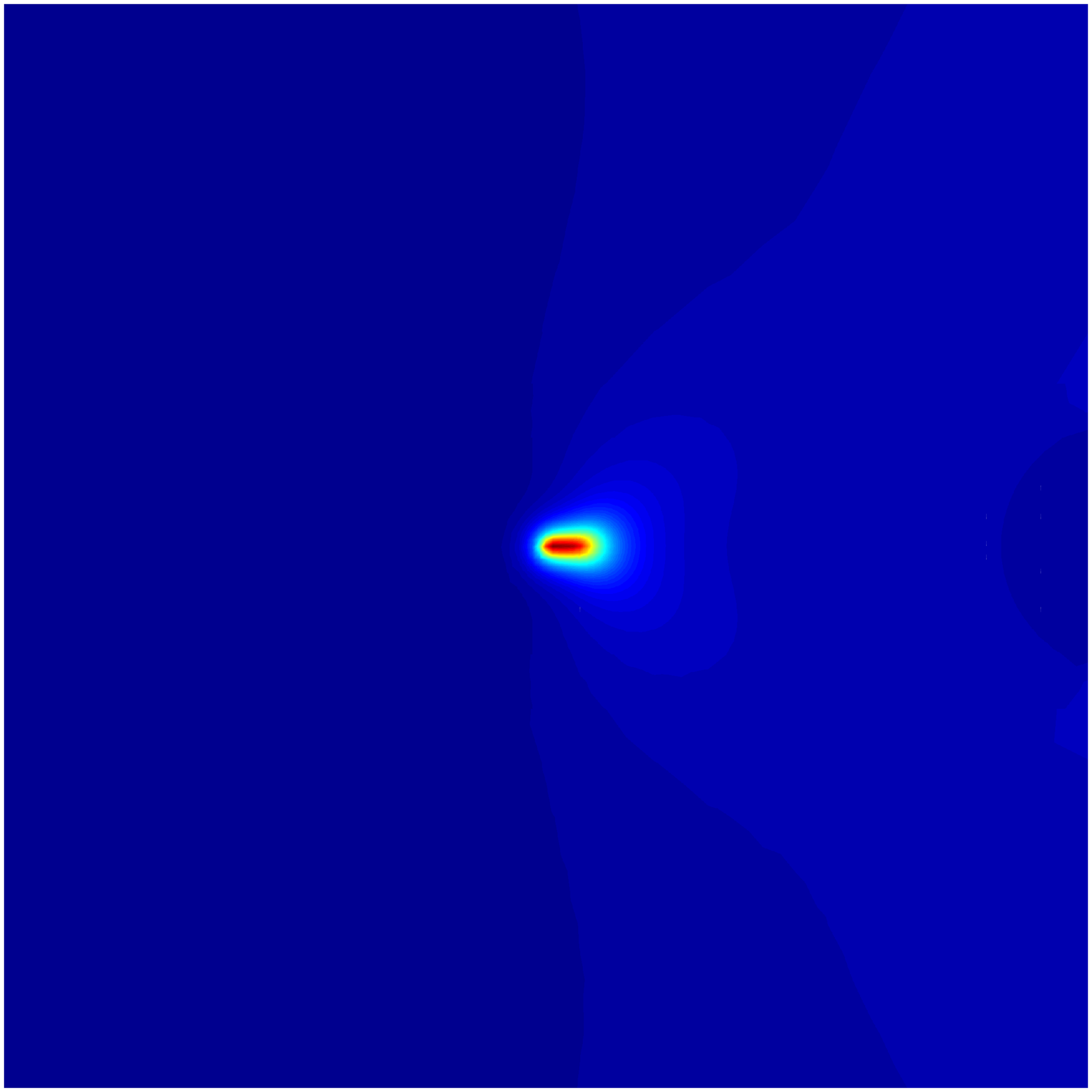}
 \includegraphics[width=0.3\textwidth]{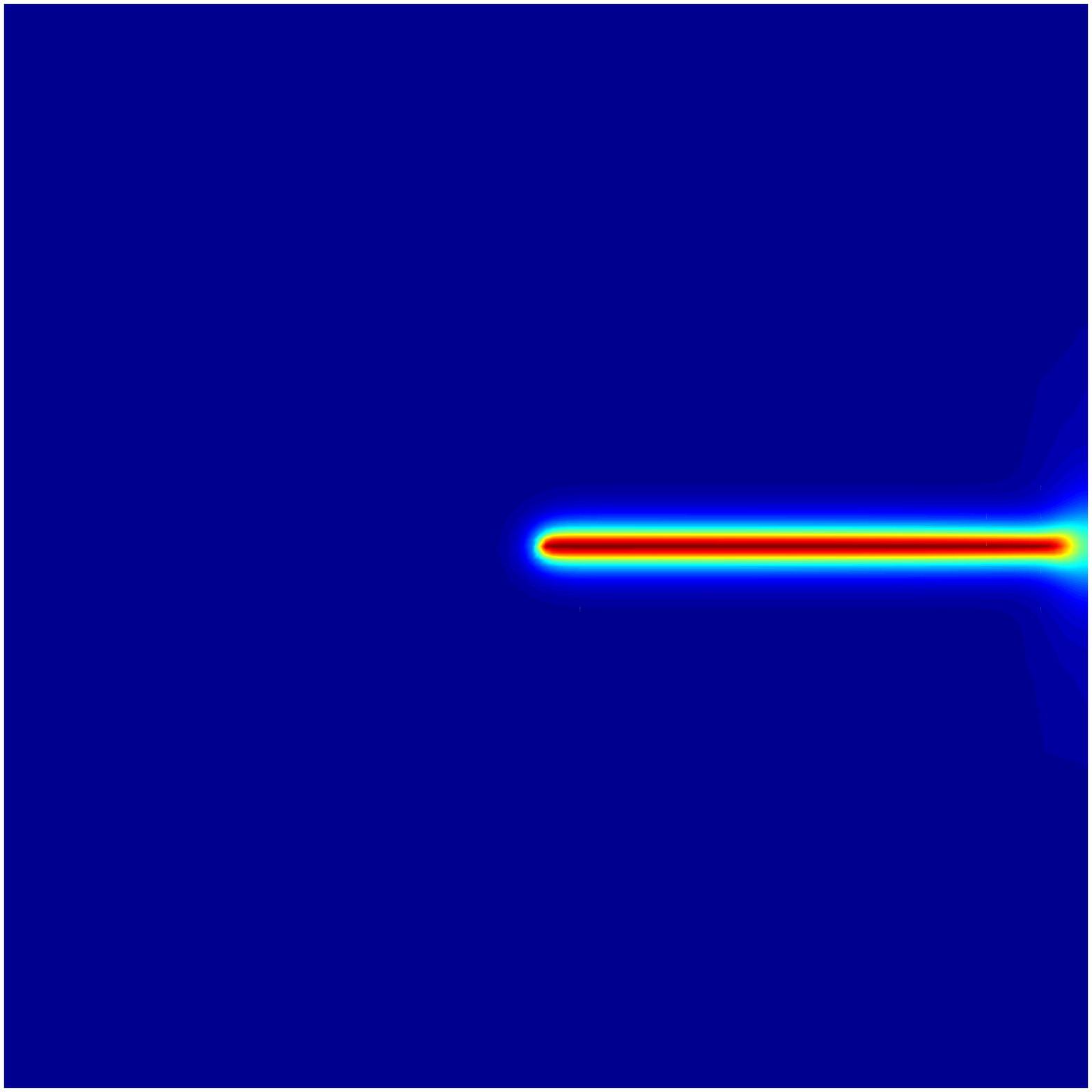}
 \caption{Geometry and boundary conditions (left) of the mode-I tension test and the typical phase-field evolution after 925 and 984 uniform steps of loading (right)}
 \label{fig:modeI}
\end{figure}

At first we  study a simple mode-I tension test and consider a $100\times100\,$mm plate with a centered notch under plane stress condition; details of geometry and   boundary conditions  are shown in Fig.~\ref{fig:modeI}.
If not mentioned otherwise, the plate is meshed uniformly with $100\times100\times4$ triangular finite elements with linear shape functions (P1 elements).   Reference material data are an elastic modulus of $ {E}=50\,400$N/mm$^2$,   Poisson ratio $\nu=0.2$ and a critical Griffith energy of $\mathcal{G}_c = 75 \, {\text{N}}/{\text{m} }$; the critical length is  $l_c=1\,$mm.
From the one-dimensional crack solution (see Appendix A) the critical stress follows as
\begin{align}\label{epsilonc}
    \sigma_c = \sqrt{\frac{E \mathcal{G}_c}{3   l_c}}
\end{align}
which gives here $\sigma_c=35\,$N/mm$^2$.
Vertical displacements   are prescribed  on the upper boundary with increments $\Delta \bar{u} = 5\cdot 10^{-5}$mm.

\subsubsection{Mode-I tension for different driving forces in linear elasticity}
We start with a comparison of the energy driven crack growth, according to the Griffith and the $\lambda$-$\mu$ split  model, and a stress driven crack growth according to the Rankine model. 
For all models the crack path is the same with the typical phase-field evolution shown in Fig.~\ref{fig:modeI}.

\begin{figure}
\centering
\psfrag{D}[cc]{\small{$\bar u$ [mm]}}
\psfrag{F}[cc][cr]{\small{$F/F_\text{max}(\lambda$-$\mu\,\text{split})$ }}
\includegraphics[width=0.45\textwidth]{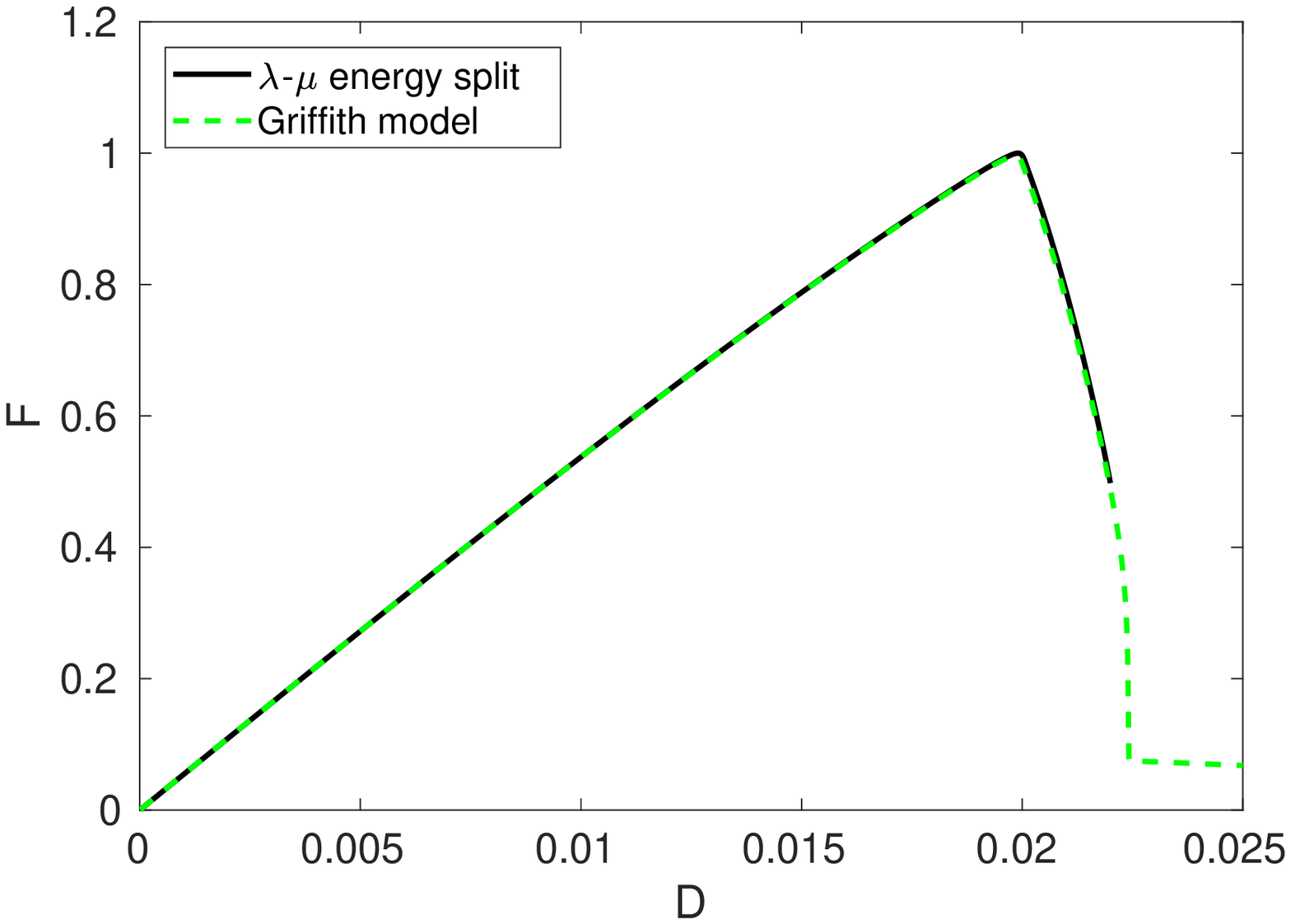}
\includegraphics[width=0.45\textwidth]{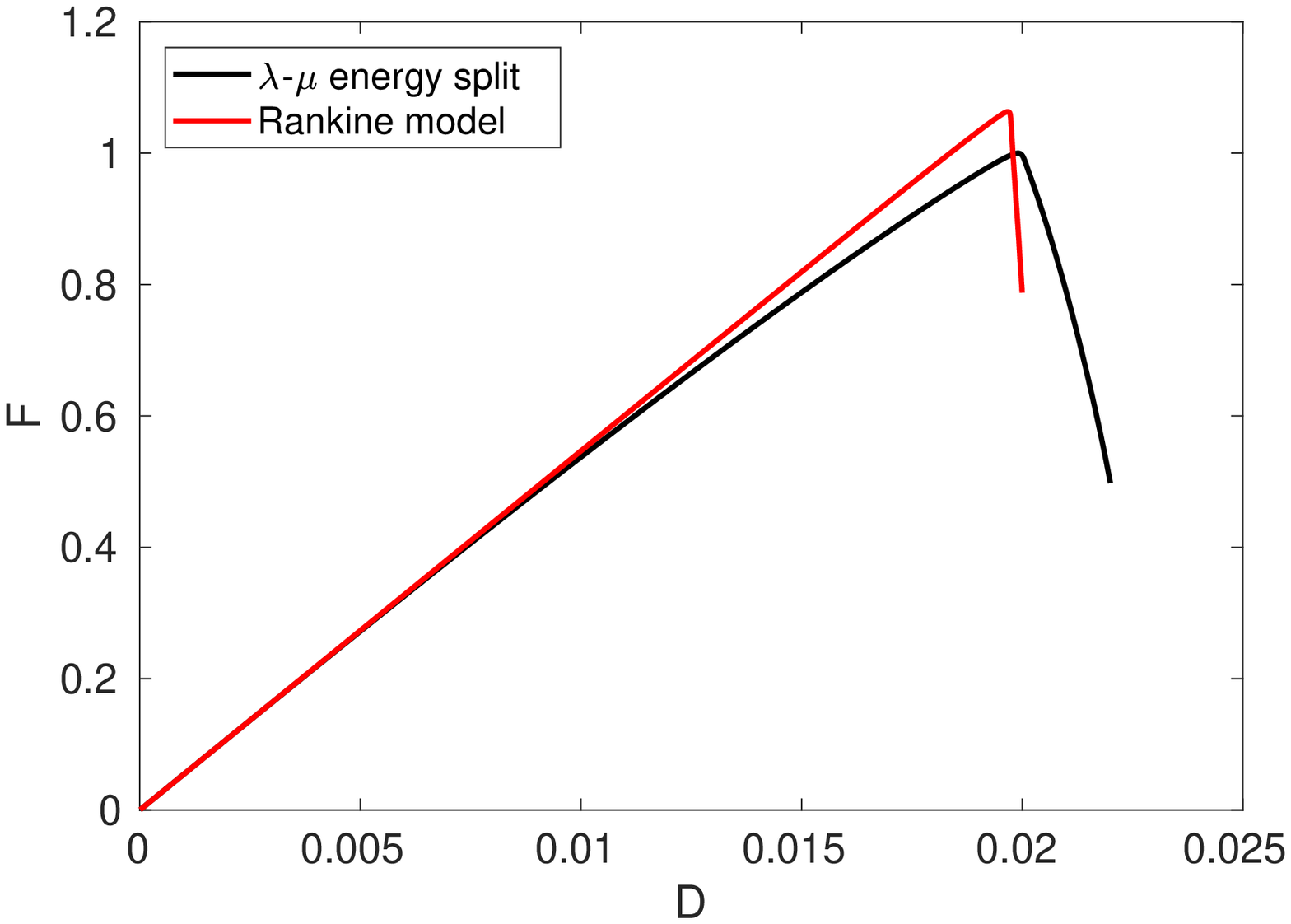}
\caption{Load-displacement curve of the mode-I tension test for two energy based driving forces (left) and a comparison of the $\lambda$-$\mu$ energy split with the Rankine model (right).  }
\label{fig:linear_LVK_diff}
\end{figure}

From the load-displacement curve in Fig.~\ref{fig:linear_LVK_diff} (left) it can be seen, that the crack starts propagating at a prescribed displacement of about $0.02\,$mm for both, the Griffith and $\lambda$-$\mu$ split  model. Both curves almost coincide which can be explained by the overall tensile state in mode I.
In Fig.~\ref{fig:linear_LVK_diff} (right) the Rankine model is compared with the energy based approaches.
They also show a very similar behavior and the crack starts propagating at a prescribed displacement of about $0.02$ mm. The remaining difference in maximum loading is caused by the definition of the critical stress via relation \reff{epsilonc}. 
Here the introduction of an additional weight for the Rankine model, as suggested in \cite{miehe2015phase}, would allow to compute identical load-displacement curves for both, the energy and the stress driven crack. This only requires a careful calibration of the additional parameter, which, however, has no physical meaningful explanation.

\subsubsection{Mode-I tension for different driving forces in  finite elasticity}
In a second step the tension test is computed in a finite deformation regime. We start with  different energy based criteria, i.e. the multiplicative split of the eigenvalues \eqref{lambdaplusoderminus} and the additive split of the invariants \eqref{SplitInv} are compared. Here we make use of a  Neo-Hookean material model extended  to the volumetric range with energy density (\ref{PsiAdditiveSplitDevVol}-\ref{MooneyRivlinNum}) with $k=0$. 
The  material data are the same as above; they result in  $K=28\,000$N/mm$^2$ and $\mu=21\,000$N/mm$^2$.

\begin{figure}
\centering
\psfrag{D}[cc]{\small{$\bar u$ [mm]}}
\psfrag{F}[cc][cr]{\small{$F/F_\text{max}(\text{energy (inv.)})$ }}
\includegraphics[width=0.45\textwidth]{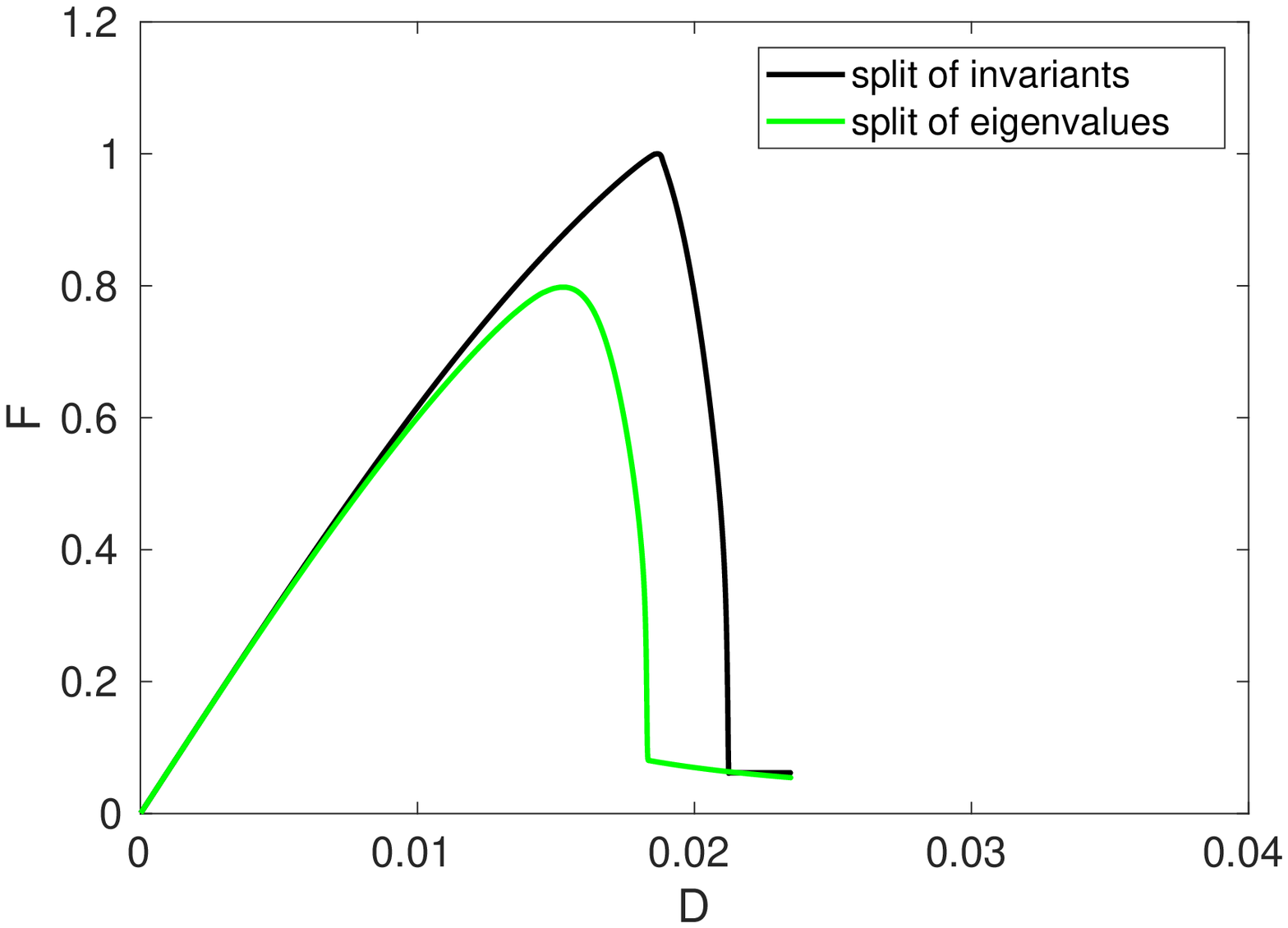}
\includegraphics[width=0.45\textwidth]{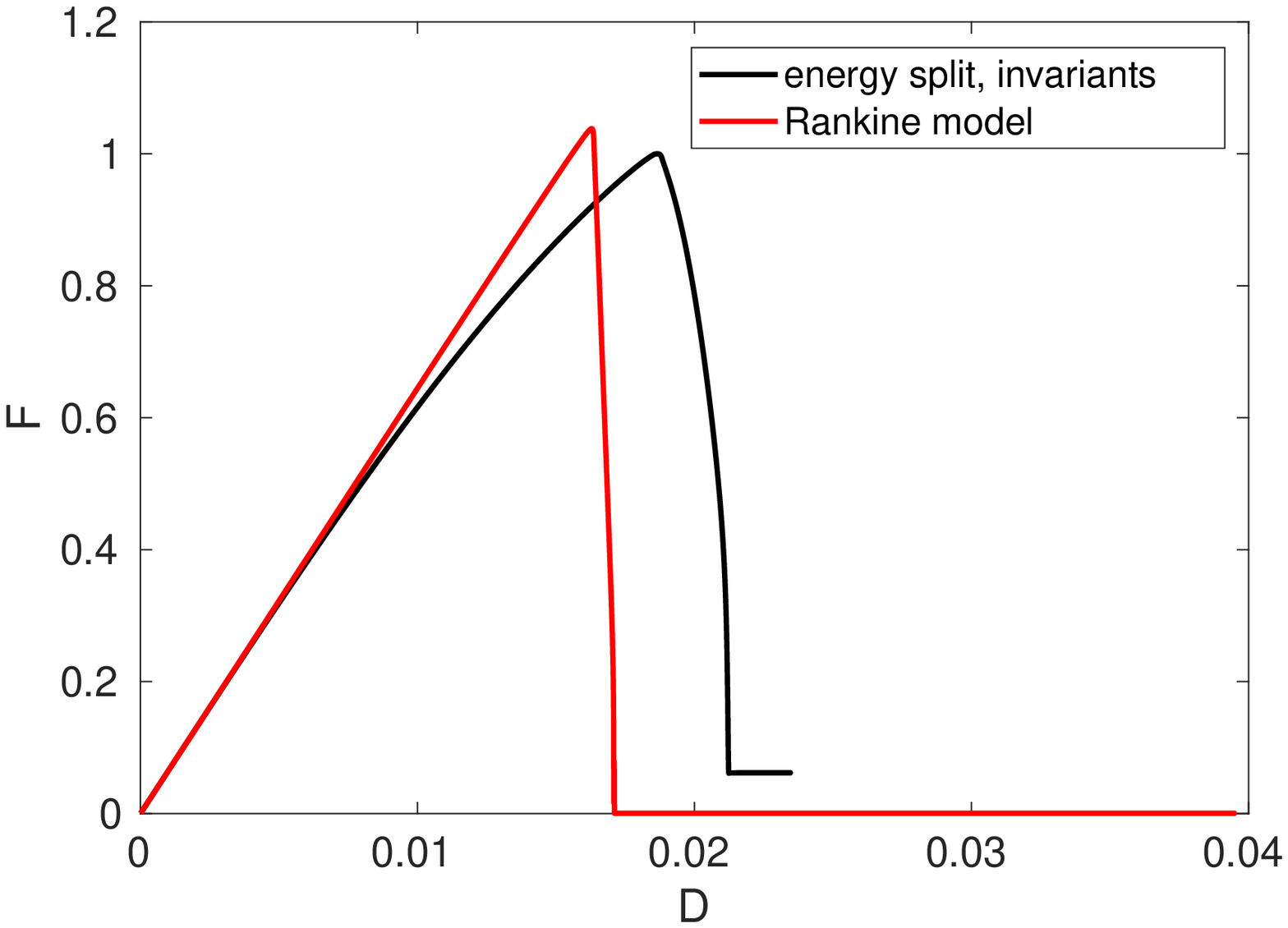}
\caption{Load-displacement curve of the mode-I tension test for different anisotropic splits (left) and for the stress based criterion (right). }
\label{fig:nonlinear_LVKdiff}
\end{figure}
The phase-field simulation of the crack path is typical mode I. The related load-displacement curves are shown in Fig.~\ref{fig:nonlinear_LVKdiff}. It can be seen that for the split of the invariants the crack starts propagating at about the same prescribed displacement as in the linear case. For the multiplicative split of the eigenvalues the load at crack initiation is about 20\% lower. Although it runs numerically stable here, the eigenvalue split seems not to be the proper choice for the crack-driving force.

The load-displacement curves of the two different driving forces -- energy split of the invariants and the Rankine model -- are depicted in Fig.~\ref{fig:nonlinear_LVKdiff} (right). Here the behavior before crack initiation is linear in the stress-based approach whereas the energy based crack-driving force is sub-linear.  For the Rankine model the crack starts propagating slightly earlier. The latter can be explained by the   critical stress value $\sigma_c$ which is derived from linear elasticity and compared here with the Cauchy stress.

\subsubsection{Effect of length scale parameter $l_c$}
At next we study the influence of  the regularization length $l_c$  on the fracture load. In the classical phase-field approach parameter $l_c$ has a two-fold effect. On the one hand, $l_c$ determines the diffusive crack width within the finite element discretization and   is   determined by the mesh size $h$. For crack evaluation $l_c>h$ is required and, depending on the approximation order, commonly  $l_c=1.5\dots 3 h$ is chosen.
On the other hand, $l_c$ enters the fracture energy density by $\ \mathcal{G}_c/l_c$. Because in the classical approach the crack will grow when the elastic energy of the material exceeds the fracture energy, $l_c$ works here like  a material parameter. In consequence, the classical phase-field approach is very sensitive against the choice of $l_c$.

\begin{figure}
 \centering
  \psfrag{force}[cc][cr]{\small{$F_\text{max}/F_\text{max}(l_c=1)$ }}
  \psfrag{displacement}[cc][cc]{\small{$\bar u$ [mm]}}
 \includegraphics[height=0.41\linewidth]{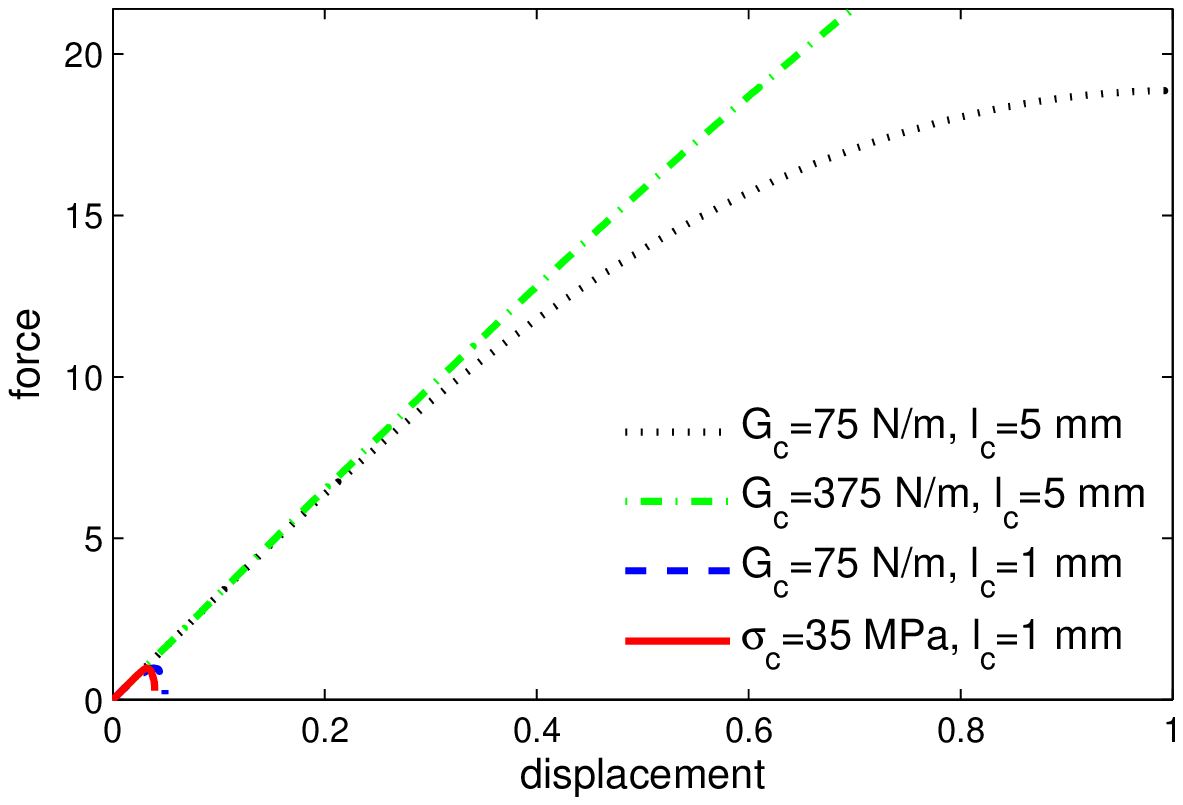}\quad
 \includegraphics[height=0.41\textwidth]{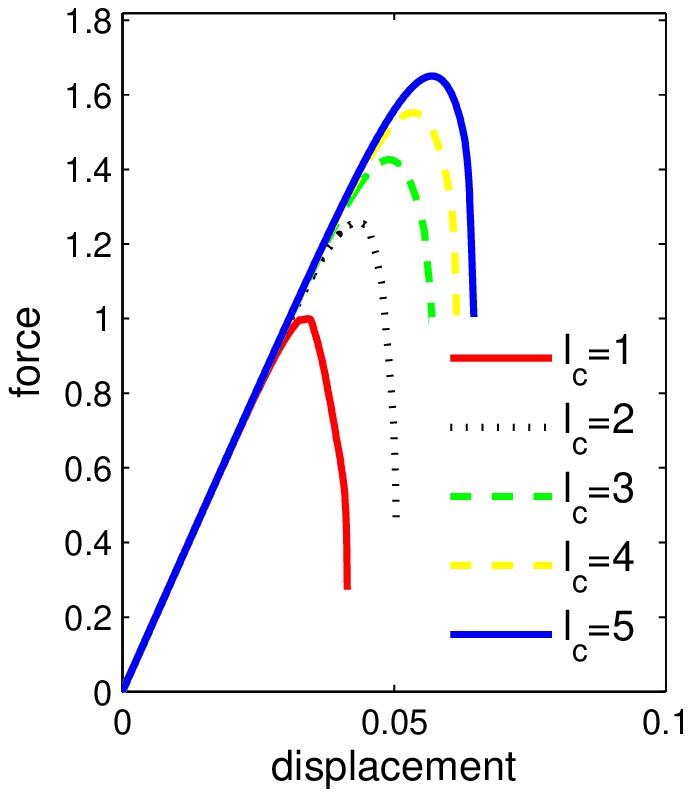}
 \caption{Variation of the critical length parameter $l_c$ for the Griffith model with constant $\mathcal{G}_c$ or $\mathcal{G}_c/l_c$ ratio (left) and for the Rankine model (right)}
 \label{fig:GcScLc1bis5}
\end{figure}

The effect can be seen in  the load-displacement curves of the linear mode-I model, plotted in Fig.~\ref{fig:GcScLc1bis5}, where we vary $l_c$ from 1 to 5\,mm.
In the Griffith model   the maximal load magnifies with $l_c$ by a factor of almost 20. This is in opposite to the effect which $l_c$ has in  stress or strain driven crack growth where $l_c$ enters only the crack width.

From the plots in Fig.~\ref{fig:GcScLc1bis5} it is obvious, that the Rankine model shows less sensitivity to  the choice of the length scale   $l_c$. In order to quantify this observation, we listed the  computed maximal loads  for different $l_c$ in Tab.~\ref{tab:cracklc} and compared it to the  diffusive 'crack volume' of a crack with unit length in the one-dimensional solution, $V(l_c)=\int_0^1 \exp(|x|/l_c) \td x$. The result is clear: the diffusive width $l_c$ effects the   load required for crack growth basically in the same way as the diffusive crack volume changes.
This is a major advantage of the Rankine crack-driving force, compared to the original phase-field model of Griffith energy driven crack growth.

\begin{table}
\begin{center}
\caption{Diffusive crack volume and maximum load  in the Rankine model for different $l_c$}\label{tab:cracklc}
\renewcommand{%
\arraystretch}{1.2}%
\begin{tabular}{|l r|ccccc|}\hline
  $l_c$ & [mm] & 1 & 2 & 3 & 4 & 5 \\  \hline 
  $V(l_c)$ 
  &[mm$^3]$ & 126 & 157 & 170 & 177 & 181 \\
  $V(l_c)/V(1)- 1 \ $ &  $ [\%]$ & 0 &  25 &  35  &  40 &  44 \\
  $F_\text{max}(l_c)\ $ &  [kN] & 93.5 &117.5  &133.3  &145.1 &  154.2 \\
  $F_\text{max}(l_c)/F_\text{max}(1) - 1\ $ & $ [\%]$ & 0 &  26 &  43 &  55&  65 \\ \hline
\end{tabular}
\end{center}
\end{table}

When in the Griffith model  the ratio $\mathcal{G}_c/l_c$ is kept constant 
the computed maximal load is even higher. This was to expect and is not displayed extra but it also underlines that the Rankine model is much more stable to modifications of $l_c$. We remark that smaller time steps would result in spikier curves. The general result, however, is not affected by a varying mesh or time discretization and it holds for all stress or strain driven crack-driving forces of Table~\ref{tab:CrackDrivingForces}.

\subsubsection{Effect of finite element discretization}
\begin{figure}
 \centering
  \psfrag{force}[cc][cr]{\small{$F_\text{max}/F_\text{max}(l_c=1)$ }}
  \psfrag{displacement}[cc][cc]{\small{$\bar u$ [mm]}}
 \includegraphics[height= 0.43\linewidth]{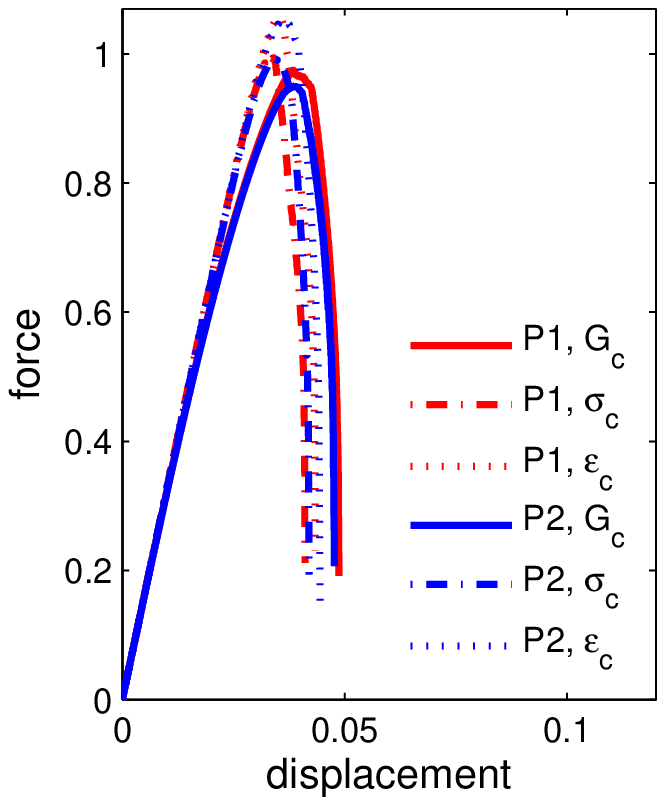}\quad
 \includegraphics[height=0.43\textwidth]{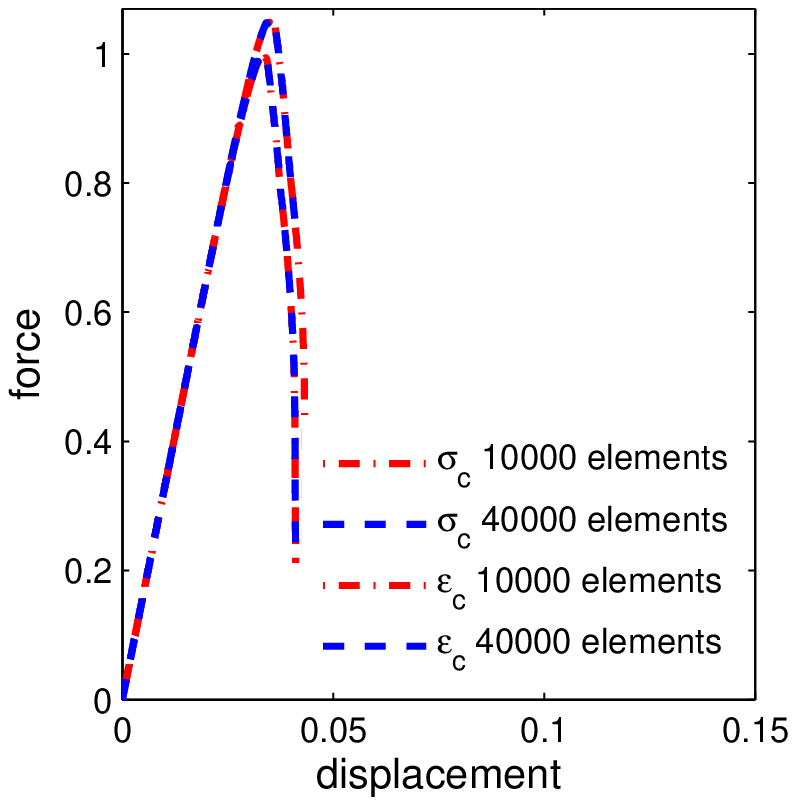}
 \caption{Effect of finite element approximation order: full energy driven, principal stress driven and principal strain driven crack growth for linear P1 and quadratic P2  triangulation (left) and principal stress  driven crack growth for different mesh sizes (right)}
 \label{fig:GcScEcP1P2}
\end{figure}

For completeness we compare now the load-displacement curves for different finite element meshes and different crack-driving forces, Fig.~\ref{fig:GcScEcP1P2}. For a fixed $l_c$ --which is   determined by the  mesh size-- we see only minor differences between a linear P1 and a quadratic P2 triangulation. A drawback of the quadratic approximation is that it tends to oscillate during unloading, affecting all phase-field models in a similar manner.
A finer mesh results in a softer response of the model and, thus, in a smaller maximum loading.
Additionally it can be seen from Fig.~\ref{fig:GcScEcP1P2} that the Rankine model and the Beltrami model give very similar results.

\subsubsection{Mode-II test}\label{sec:modeItensionTest}
\begin{figure} 
 \centering
  \psfrag{U}[cc]{$\bar{u}$}
  \psfrag{w}[cc][cc]{$\alpha$}
  \psfrag{A}[cc][bc]{Griffith model}
  \psfrag{B}[cc][cc]{principal axis   split}
  \psfrag{C}[cc][cc]{$K$-$\mu$ energy split}
  \psfrag{D}[cc][cc]{$\lambda$-$\mu$ energy split}
  \psfrag{E}[cc][cc]{Tresca model}
  \psfrag{F}[cc][cc]{Rankine model}
  \psfrag{G}[cc][cc]{compressive Rankine }
  \psfrag{H}[cc][cc]{Mohr-Coloumb model}
  \psfrag{J}[cc][cc]{Beltrami model}
  \psfrag{K}[cc][cc]{inverse Mohr-Coloumb model}
  \includegraphics[width=0.99\linewidth]{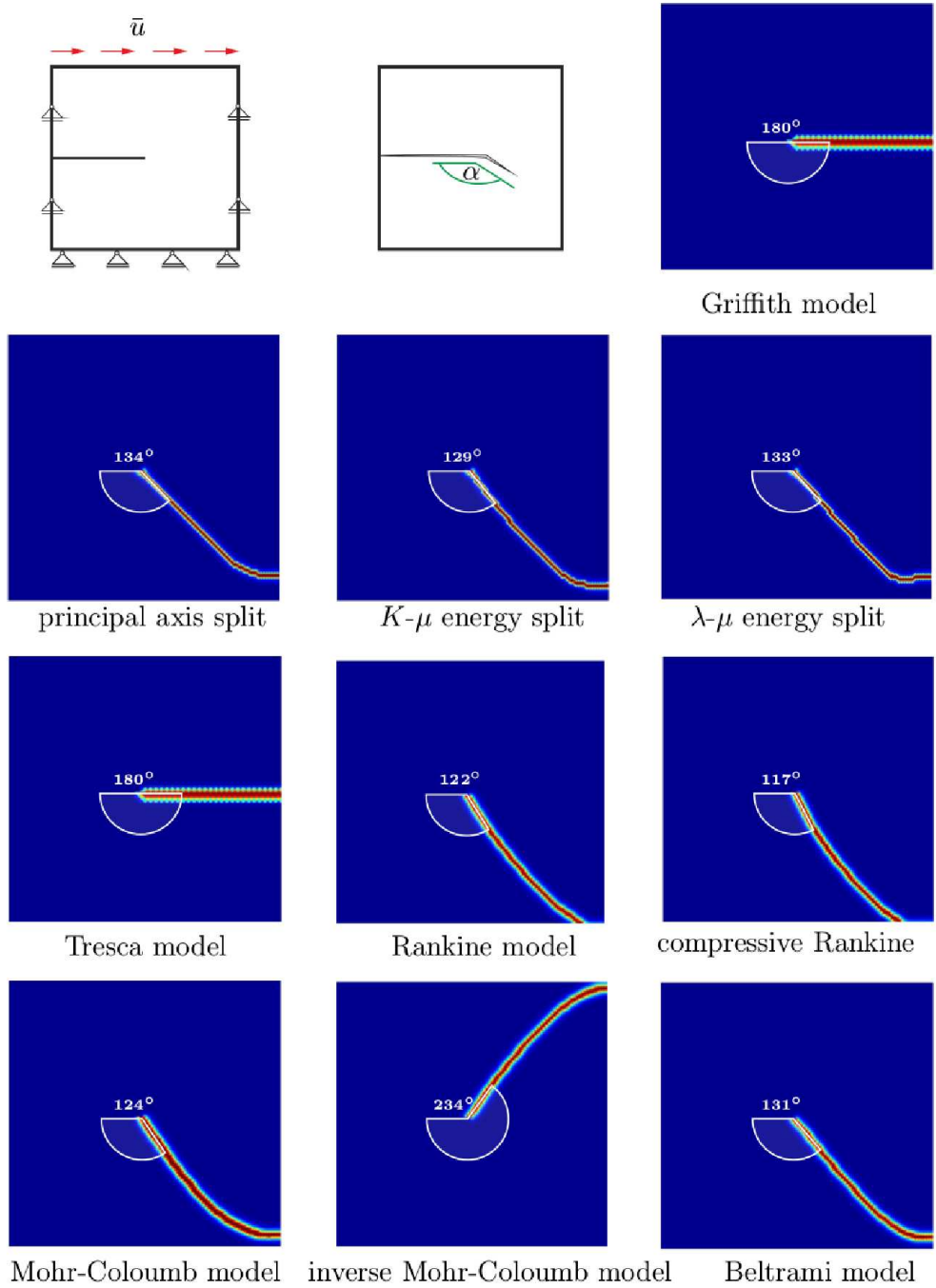}
 \caption{Mode-II-shear test with different crack-driving forces}
 \label{fig:modeII12}
\end{figure}

Because the tension test is  too simple to compare the different crack-driving forces, we change the boundary conditions of our computational model now to a shear test.
This situation has been used for numerous investigations on the direction of kinking cracks. 
Different theoretical approaches have been developed, e.g. the principle of local symmetry, which states that the crack propagates in such a way that the mode-II stress intensity factor vanishes in the vicinity of the crack tip, \cite{Amestoy1992,Goldstein1974,Hodgdon1993,Leblond1989,katzav2007fracture}.  
Alternatively the direction of crack growth can be determined so that it minimizes the potential energy or maximizes the energy release rate among all kinking angles, \cite{Erdogan1963,Bilby1975,Chambolle2009}. Experiments under shear have been performed, for example, by Erdogan and Sih \cite{Erdogan1963}. They  observed  a crack declined by an angle of $\alpha=110^\circ$ in sheared PMMA 
and mentioned also, that the kinking direction strongly depends on the 'brittleness' of the material.

In this sense we study here the effect of all crack-driving forces in shearing. We start with a linear elastic plane stress computation and observe   different kinking angles, Fig.~\ref{fig:modeII12}.
The full elastic energy of the Griffith model would drive a straight crack growth which is shear dominated and does not kink.
The different energy splits give similar and realistic results, but in detail the kinking angles differ between $129^\circ$ and $134^\circ$. The Tresca model, which accounts for shear only, gives a straight $180^\circ$ crack growth. This may be observed in ductile metals. The Rankine model of maximal principal stresses  results in a kinking angle of $122^\circ$, the compressive Rankine model 
reduces this angle to $117^\circ$ and gives the steepest kink. Quite interesting results gives the Mohr-Coloumb model which accounts for tension and shear. If we assume a tension sensitive material ($m=10$) the crack kinks as in the previous models with $124^\circ$ whereas a compression sensitive material ($m=0.1$) like, e.g. foam, shows a completely opposite behavior with a shear induced crack towards the upper side of the model. At last, the  Beltrami model of maximal principal strains behaves similar to the principal stress model with a kinking angle of $131^\circ$.

\begin{figure}
 \centering
 \includegraphics[width=0.85\textwidth]{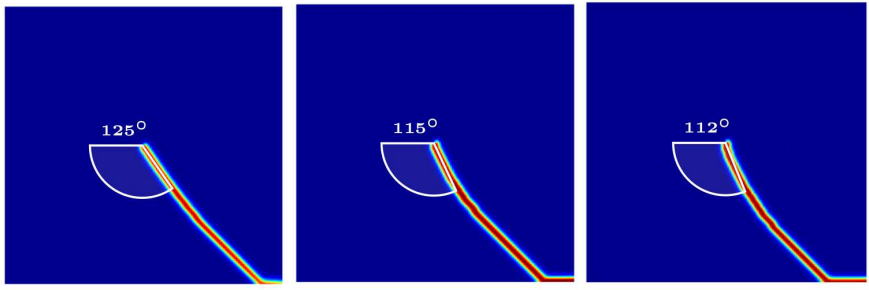}
 \caption{Phase-field snapshots of the two dimensional sheartest for the different models - anisotropic split of the eigenvalues with $\alpha = 125°$ (left), anisotropic split of the invariants with $\alpha = 115°$ (middle) and the Rankine-model with $\alpha=112°$ (right).}
 \label{fig:nonlinear_ShearAngle}
\end{figure}

In Fig.~\ref{fig:nonlinear_ShearAngle} the crack evolution in the mode-II shear test for the different finite elasticity models is displayed. Here again, the different energy splits give different results but the best result ---in comparison to the experiments in \cite{Erdogan1963}--- gives the Rankine model of maximum principal stresses. For all models the stress driven crack growth is similar to the linear elastic regime whereby the kinking angle is always somewhat smaller.   This may in part be due to the fully three-dimensional computation.

In total we may summarize, that the generalized driving force model offers the opportunity to describe crack characteristics, like as the kinking angle under shear, as it is observed in practise. The right choice of the crack-driving force depends on the  properties of the specific material.

\subsection{Brazilian test simulations}
The Brazilian  or compressive-split test of civil engineering is a popular experiment to determine the tensile strength of brittle, tension-sensitive material. Two opposing forces are applied to the cross section of a cylindrical specimen until a tension, perpendicular to the loading direction, causes the specimen to crack. The Brazilian test   with quasi-static loading is standardized, e.g. in ASTM C496 and  DIN EN 12390-6, 
and determines the tensile resistance of the specimen as
\begin{equation}\label{tensilesplittingstrength}
R_m^\text{t}=\frac{2F_\text{max}}{\pi L D}
\end{equation}
where $F_\text{max}$, $L$ and $D$ are maximum applied load, length and diameter of the specimen, respectively.

\begin{figure}
 \psfrag{D}[cl][cc]{\small{$D=50\,$mm}}
 \psfrag{b}[cc][cc]{\small{$b$}}
 \psfrag{u}[cc][cc]{\small{$\bar{u}$}}
 \includegraphics[width=0.28\linewidth]{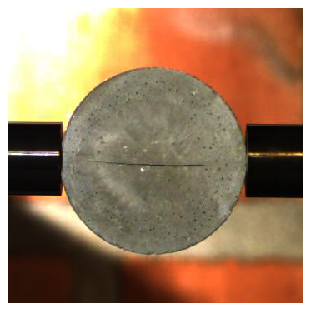}\hfil
\includegraphics[width=0.37\textwidth]{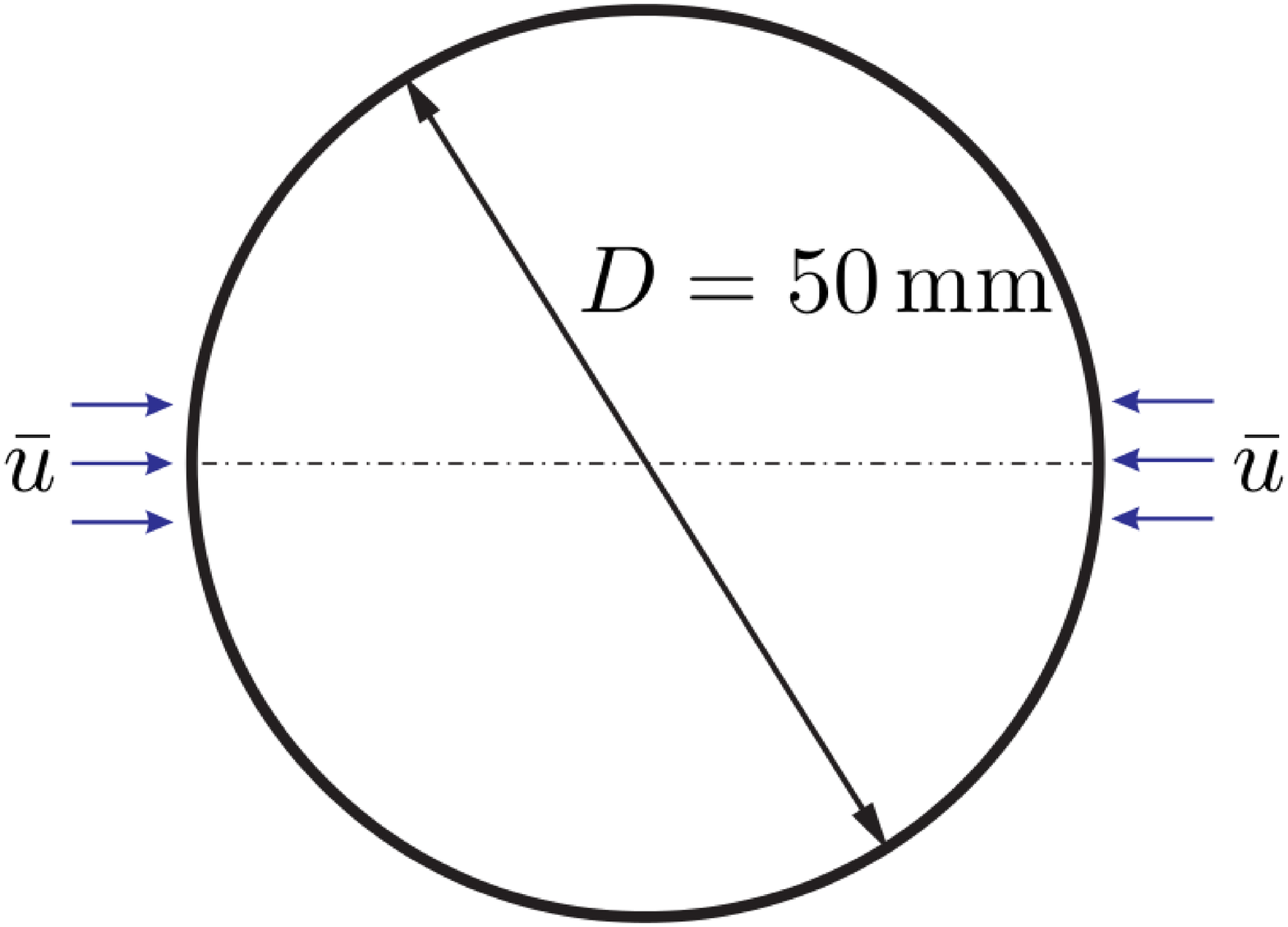}\hfil
\includegraphics[width=0.27\textwidth]{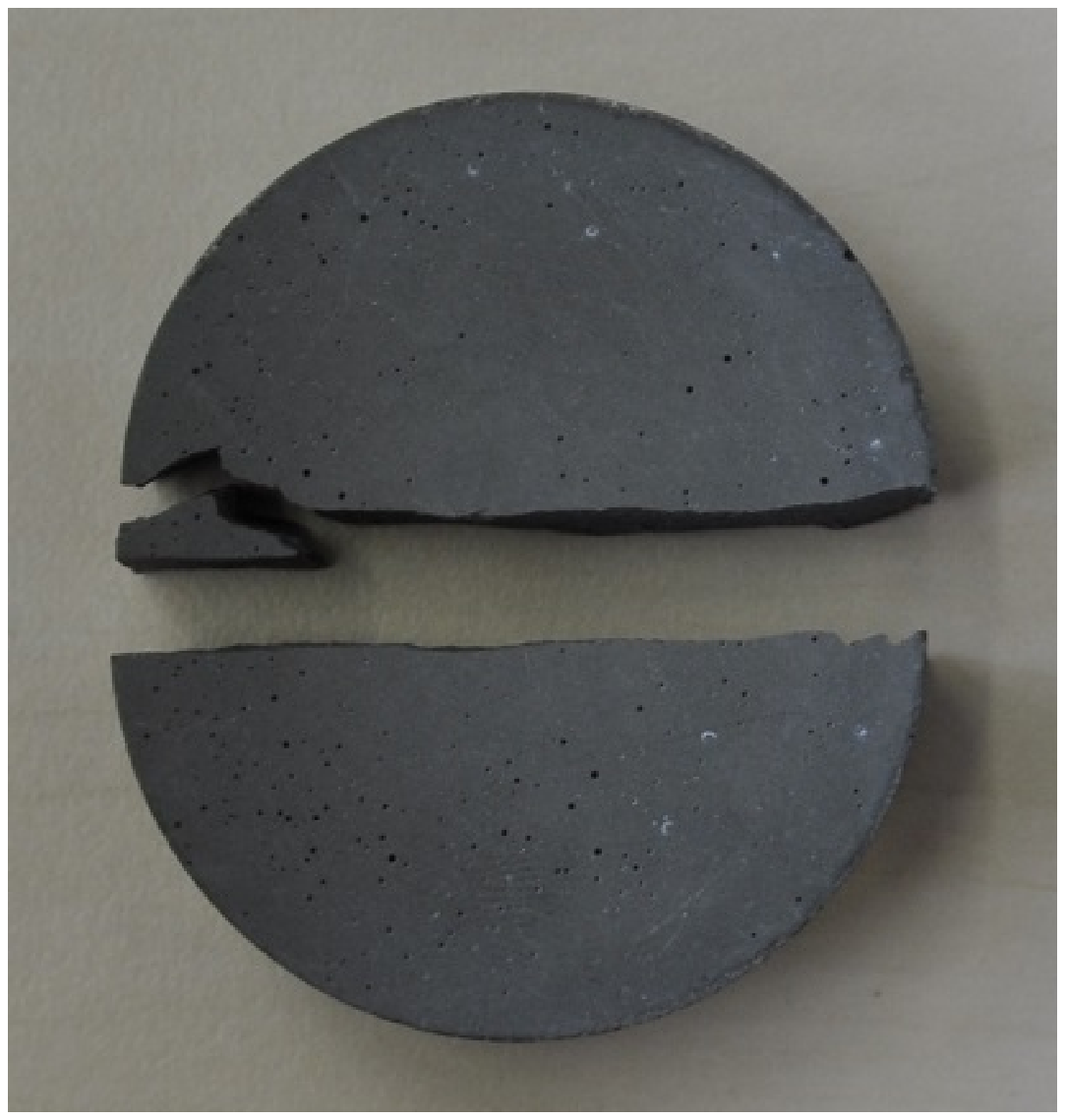}
\caption{UHPC Brazilian disc specimen: dynamic cracking (left), specimen geometry (middle) and cracked specimen geometry (right)}
\label{fig:brazspecgeometry}
\end{figure}

\begin{figure}
 \psfrag{A}[cc][cc]{\small{Griffith model}}
 \psfrag{B}[cc][cc]{\small{$\lambda$-$\mu$-split}}
 \psfrag{C}[cc][cc]{\small{Rankine model}}
\includegraphics[width=0.45\textwidth]{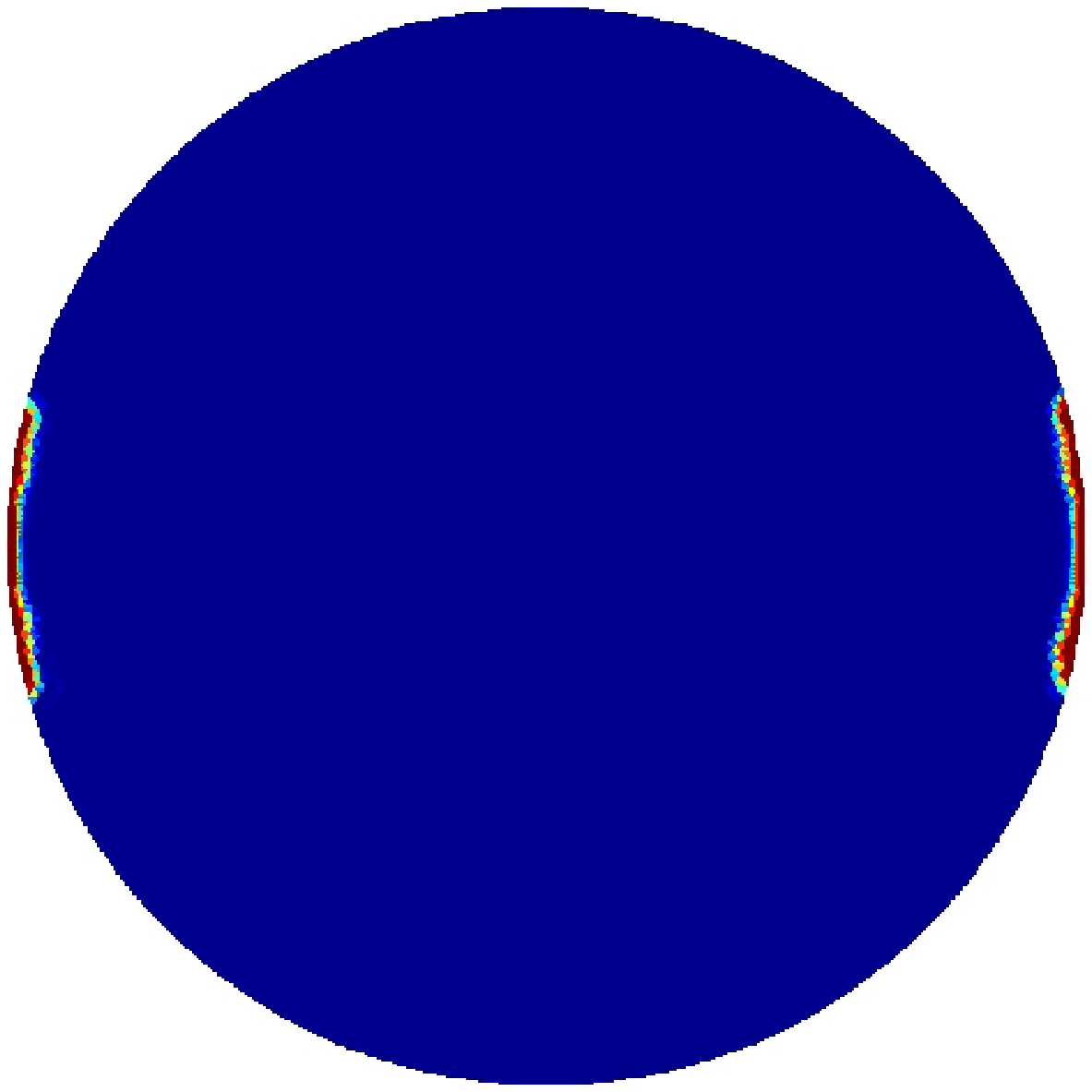}\hfil
\includegraphics[width=0.45\textwidth]{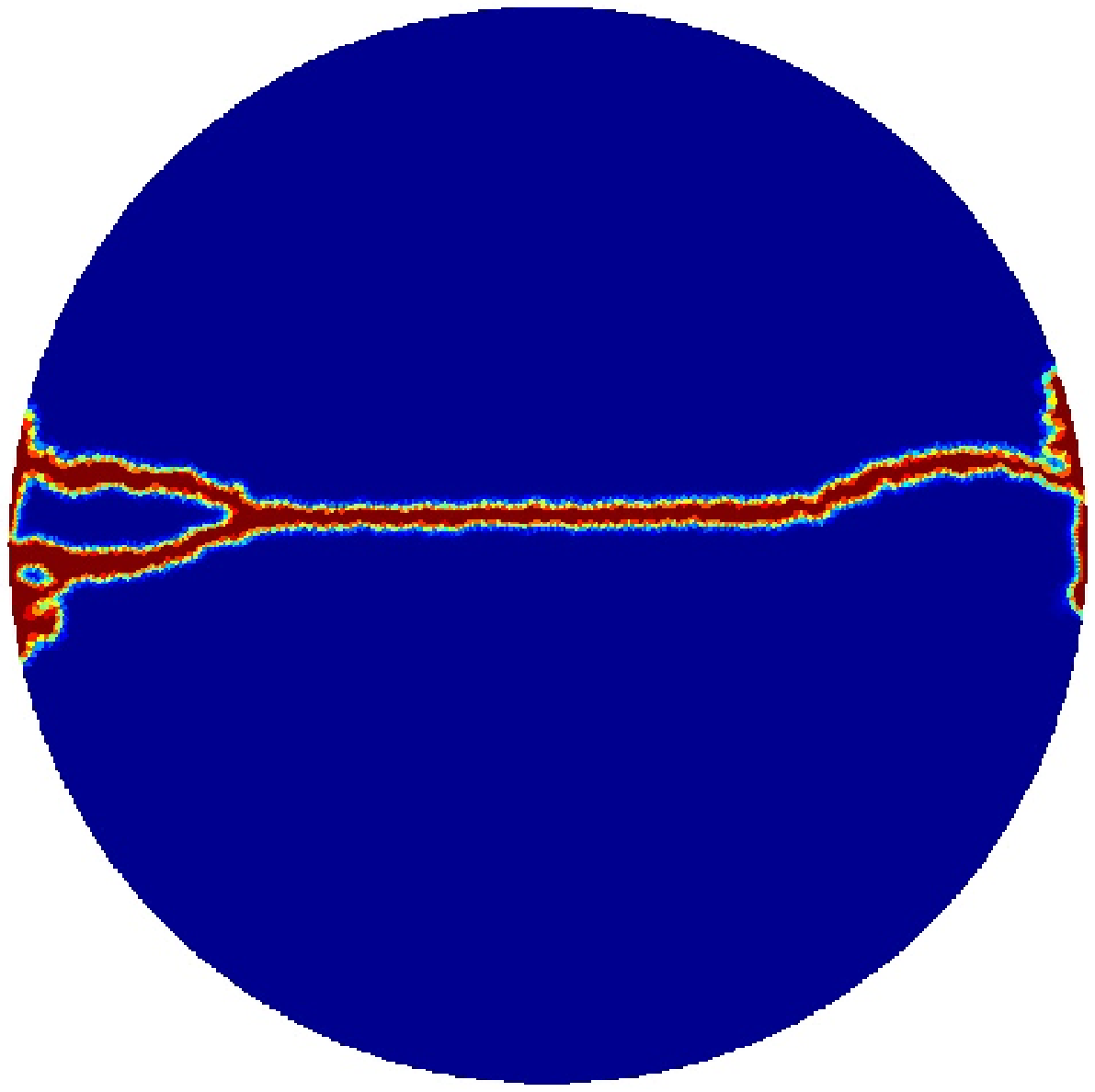}
\caption{Crack pattern of the Brazilian disc specimen for different crack-driving forces: Griffith model (left) and Rankine model (right)}
\label{fig:BrazilianTMP}
\end{figure}

In \cite{Khosravanietal2018TAFM} we performed dynamic Brazilian test with specimens made of   Ultra-High Performance Concrete (UHPC), whereby the evaluation of the  tensile resistance is the same as in statics, Eq.~\reff{tensilesplittingstrength}.  Setup,  geometry and a cracked specimen are displayed in Fig.~\ref{fig:brazspecgeometry}  and serve us here as a reference to study phase-field fracture simulations.  Our finite element model consists of {33\,169} P1 elements.  A successive compressive  displacement in axial direction is prescribed from both sides. At the two center points also the lateral displacement is constraint. The specimen has the typical material data of UHPC, $E=50\,000\,$MPa, $\nu=0.2$, $R_m^\text{t}=  \sigma_c = 17\,$MPa and $G_c=70\,$N/m.

Please note that the analytical derivation of equation \reff{tensilesplittingstrength} presumes a circular linear-elastic disc with point loads $F$,  \cite{Frocht1948}. The   principal tensile stress, $\sigma_I$, which has its maximum in the center of the disc,  is then identified with the resistance \reff{tensilesplittingstrength}. 
The value of the compressive   stress in the specimen is much higher.
From the practical point of view, the major criterium for a valid Brazilian test is that the crack starts in the center of the disc. This happens only for materials with a strong tension-compression asymmetry. Other brittle materials like glass or some plastics fail the   test because their cracking starts at the support.

Figure~\ref{fig:BrazilianTMP} shows the computed specimen for different crack-driving forces. All the energy criteria fail to give realistic crack patterns.
In our computation this means at the loading zone a rapid loss of stiffness down to the $\varepsilon$-value of the degradation function \reff{weightf}. Further results are then meaningless.
Realistic crack pattern in the sense of a valid Brazilian test gives only the Rankine model and, slightly different, the compressive Rankine model. In  both models the crack starts in the center of the disc and grows rapidly towards the boundary.

The maximum applied force per thickness in the Rankine model is  $F_\text{max}=3300\,$N/mm which is roughly in the range of the value $F_\text{max}= 2750\,$N/mm expected from Eq.~\reff{tensilesplittingstrength}. The calculated force depends on the crack-driving force but also on the boundary conditions and will be elaborated in more detail in a subsequent work. Here we summarize that, regarding the theoretical basis of the experiment, only the phase-field crack-driving force determined  by the maximum tensile stress   models the Brazilian test  adequately.


\subsection{A conchoidal fracture example}
Finally a three dimensional example is studied to demonstrate the effect of different crack-driving forces. We focus on a specific type of brittle fracture known as conchoidal fracture. Conchoidal fracture can be observed in fine-grained or amorphous materials, such as rock and glass, and  is characterized by  cleavage without any natural planes of separation. It results in a curved surface with ripples, the so called Wallner lines \cite{Wallner1939}, that resembles the surface of a sea shell.
In \cite{Bilgen_etal2017}    phase-field fracture computations are performed on  two conchoidal fracture examples and a series of investigations were carried out, mainly focussing on a fast numerical scheme for the classical approach. Here we want to extend these results by applying different crack-driving forces in the phase-field model.

\begin{figure}
\centering
\includegraphics[width=0.5\textwidth]{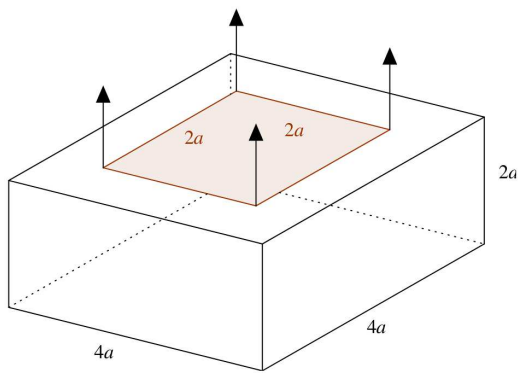}
\caption{Geometry and loading of a bloc of brittle material. }
\label{fig:BM_geometry}
\end{figure}

The geometry of the problem is illustrated in Fig.~\ref{fig:BM_geometry}. A three-dimensional $4a\times4a\times2a$ bloc of brittle   material   is loaded by a vertical displacement  on part of its upper boundary, $2a=1\,$m. On the lower boundary, Dirichlet conditions $\V u=0$ and $z=0$ for the displacements and the phase-field are prescribed.

\subsubsection{Linear-elastic analysis}
The stone-like material has an elastic modulus of $E =250000\ \nicefrac{\text{N}}{\text{mm}^2}$, a Poisson's ratio of $\nu=0.25$ and a Griffith energy of $\mathcal{G}_c = 1\ \text{N}/\text{mm}$. We use a finite element mesh which consists of $27 000$ 8-node brick elements, such that the geometry can be resolved with a length scale parameter of $l_c = 0.2\,$m.
The displacement-driven deformation with increments of $\Delta \bar{u} =0.0001\,$mm is prescribed up to complete failure.

\begin{figure}
\centering
\includegraphics[width=0.3\textwidth]{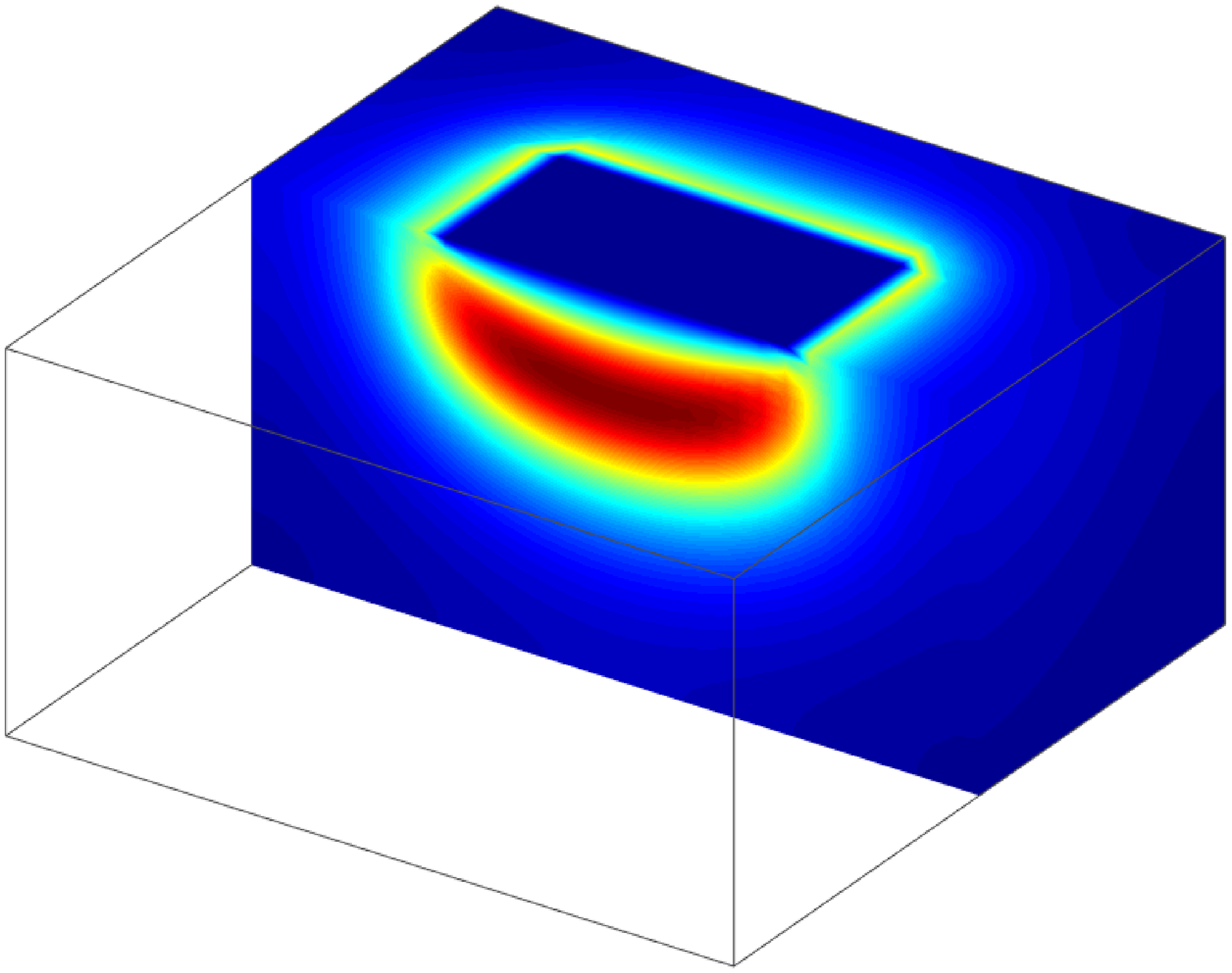}
\includegraphics[width=0.3\textwidth]{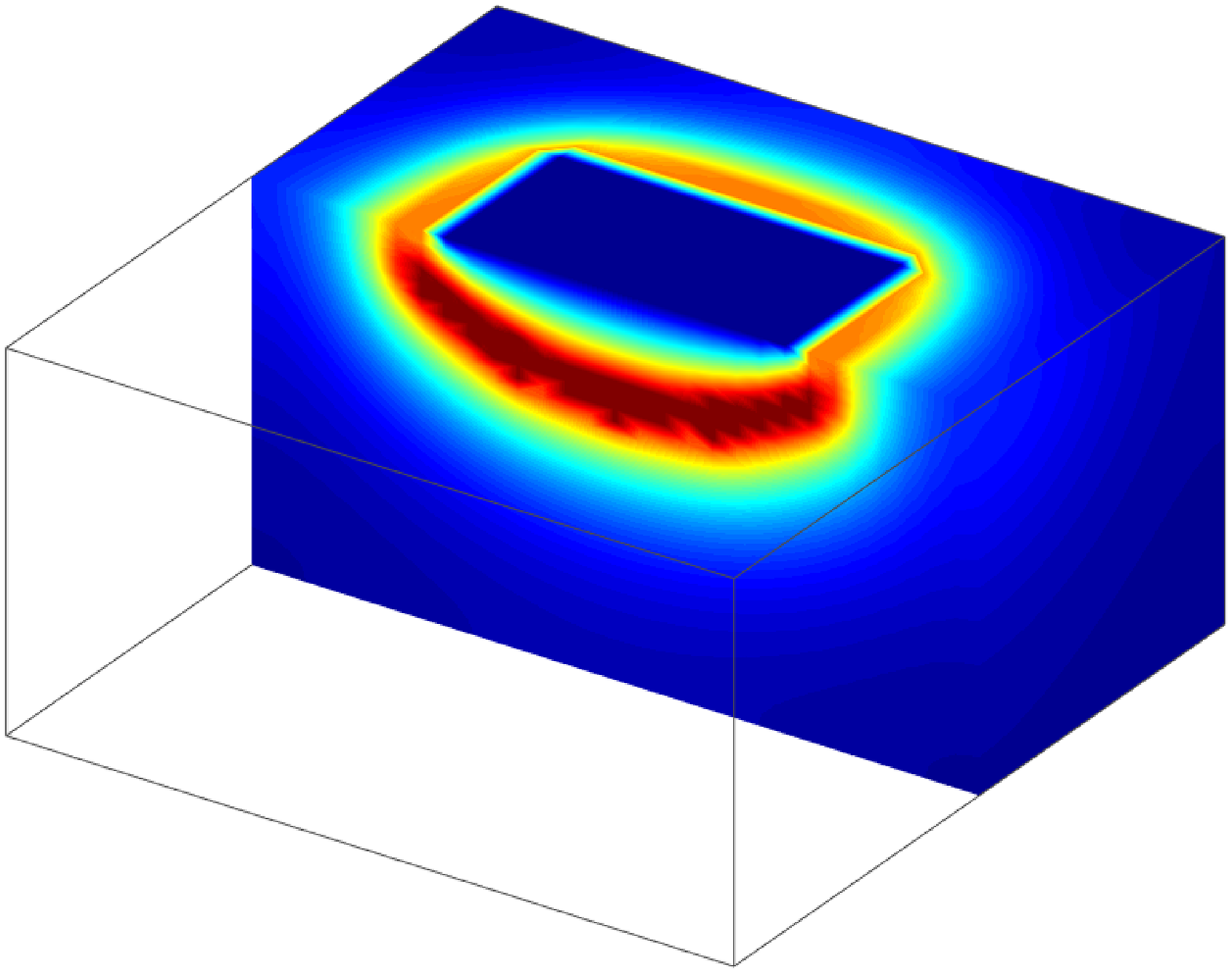}
\includegraphics[width=0.3\textwidth]{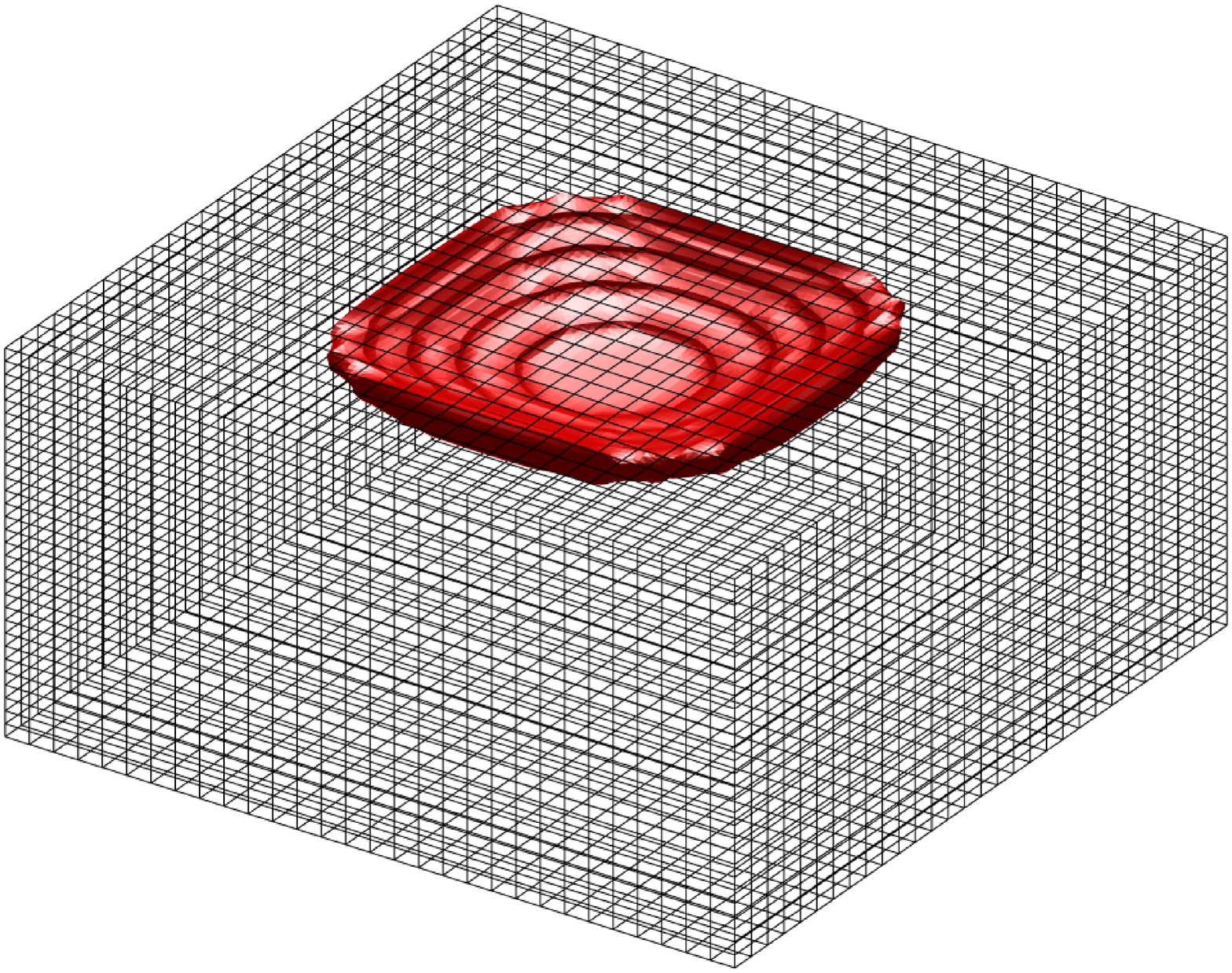}
\caption{Phase-field snapshots applying the Rankine model formulation at time step 20 (left), time step 30 (middle) and iso-surface at time step 30 (right).}
\label{fig:BM3D_stress_PF}
\end{figure}

Snapshots of the phase-field evolution at different time steps can be seen in Fig.~\ref{fig:BM3D_stress_PF}. For the  displayed maximum-stress driven crack (Rankine model) as well as for all energy-split driven models the crack starts propagating inside of the bloc below the pulled surface. Please note that this crack nucleation happens in a completely homogenous material without any notch or flaw. Once nucleated the crack is characterized   by  brutal  growth until complete failure.

\begin{figure}
\centering
\psfrag{D}[cc]{\small{$\bar u$ [mm]}}
\psfrag{F}[cc][cr]{\small{$F$ [N] }}
\includegraphics[width=0.7\textwidth]{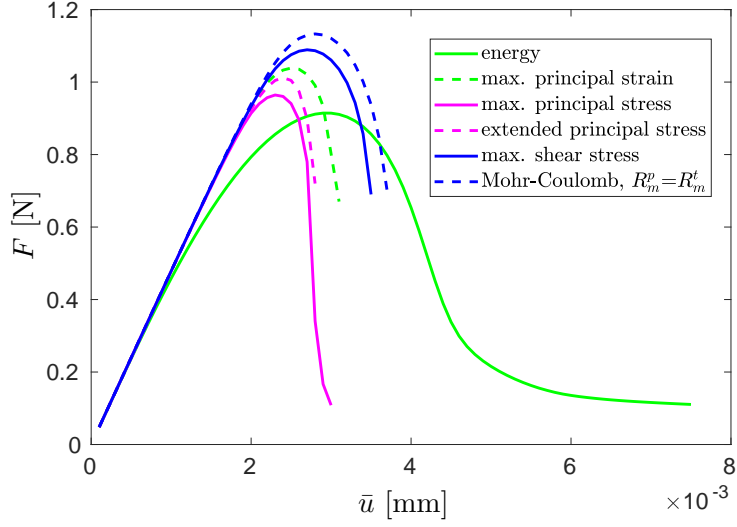}
\caption{Load-displacement curves of the conchoidal fracture example for   different crack-driving forces. The names in the legend refer to the models of  Table~\ref{tab:CrackDrivingForces}.}
\label{fig:BM3D_criteria}
\end{figure}

In the iso-surface 
the typical  rippled surface of conchoidal fracture can be observed. The crack nucleates at a displacement of approximately $0.003\,$m , see Fig.~\ref{fig:BM3D_criteria}. As already noticed in the previous examples the maximum load of crack initiation differs slightly for the different driving forces.   

%

\begin{figure}
\centering
\psfrag{D}[cc]{\small{$\bar u$ [mm]}}
\psfrag{F}[cc][cr]{\small{$F$ [N]}}
\includegraphics[width=0.7\textwidth]{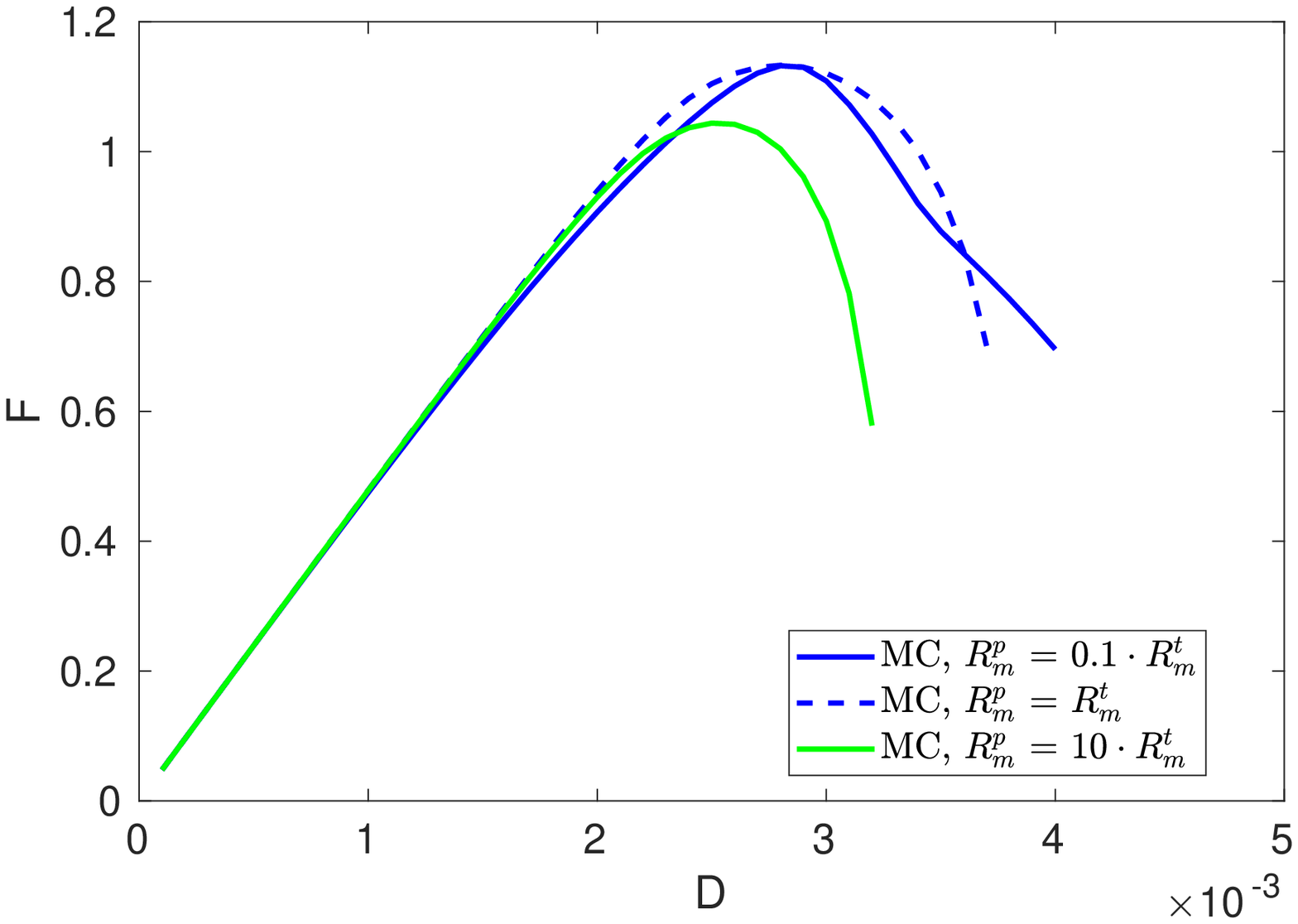}
\caption{Load-displacement curves of the conchoidal fracture example for the Mohr-Coulomb crack-driving force with different ratios of $R_m^p/R_m^t=m$.}
\label{fig:BM3D_MC}
\end{figure}

\begin{figure}
\centering
\includegraphics[width=0.3\textwidth]{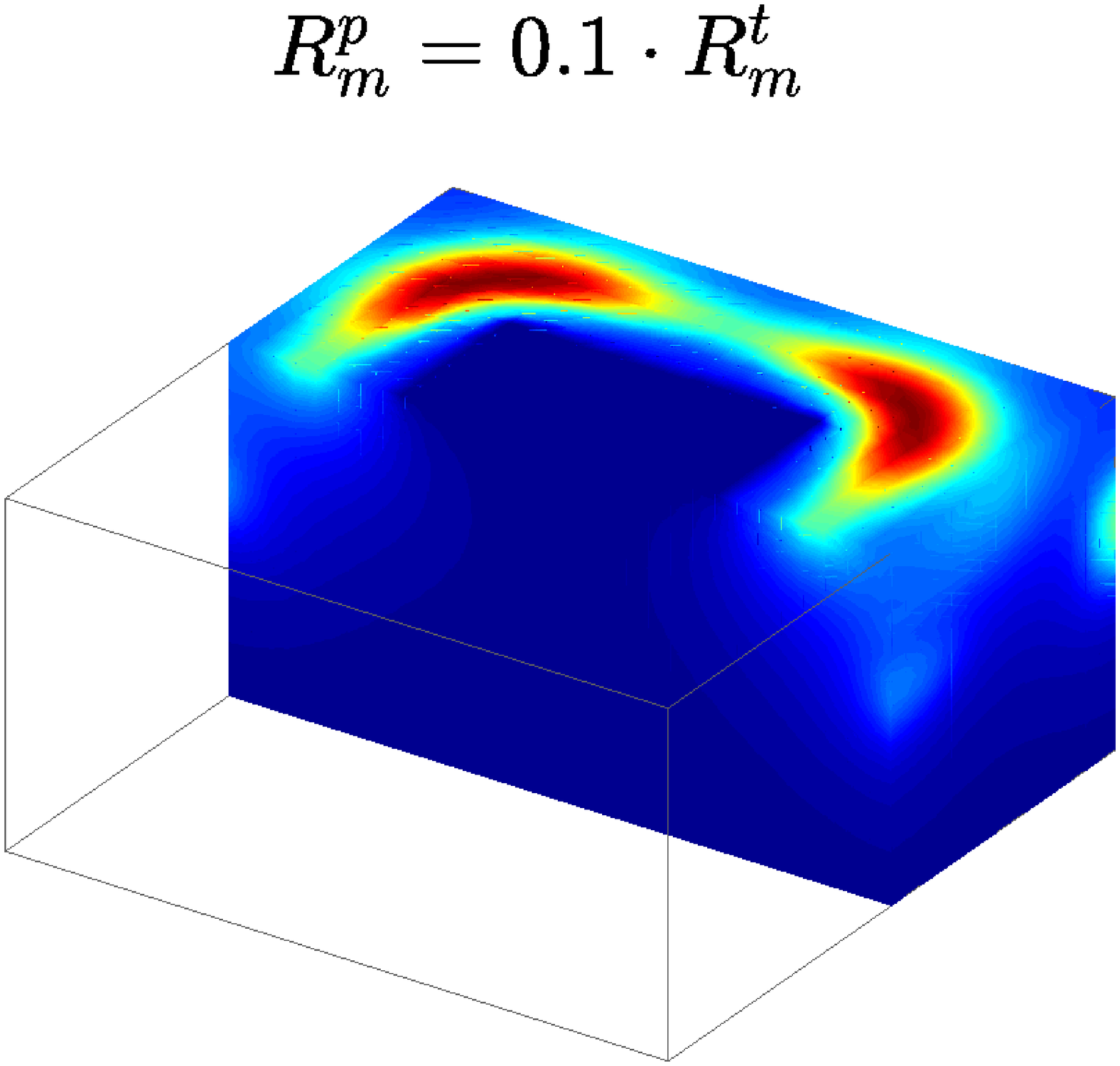}\, \includegraphics[width=0.3\textwidth]{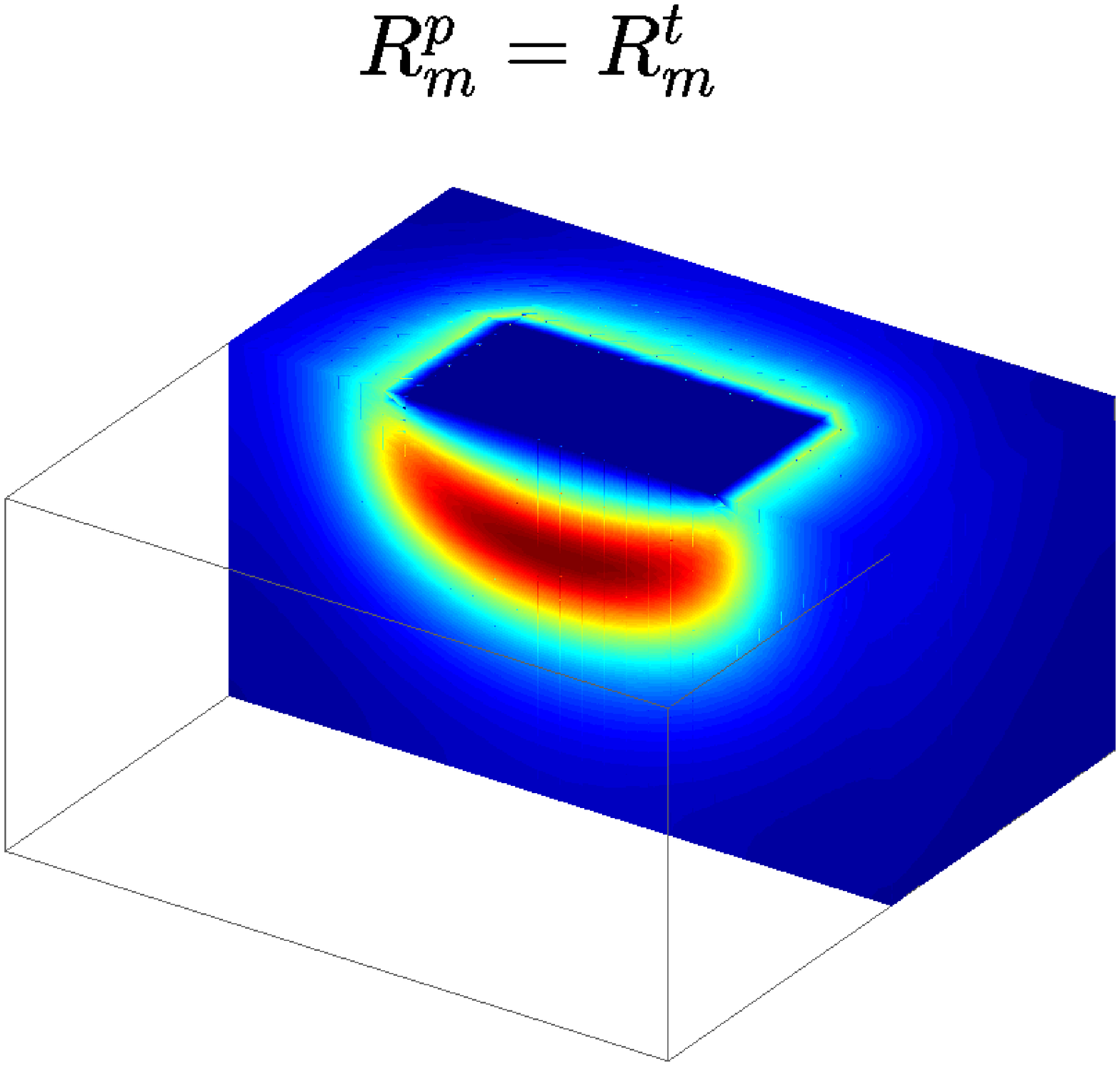}\, \includegraphics[width=0.3\textwidth]{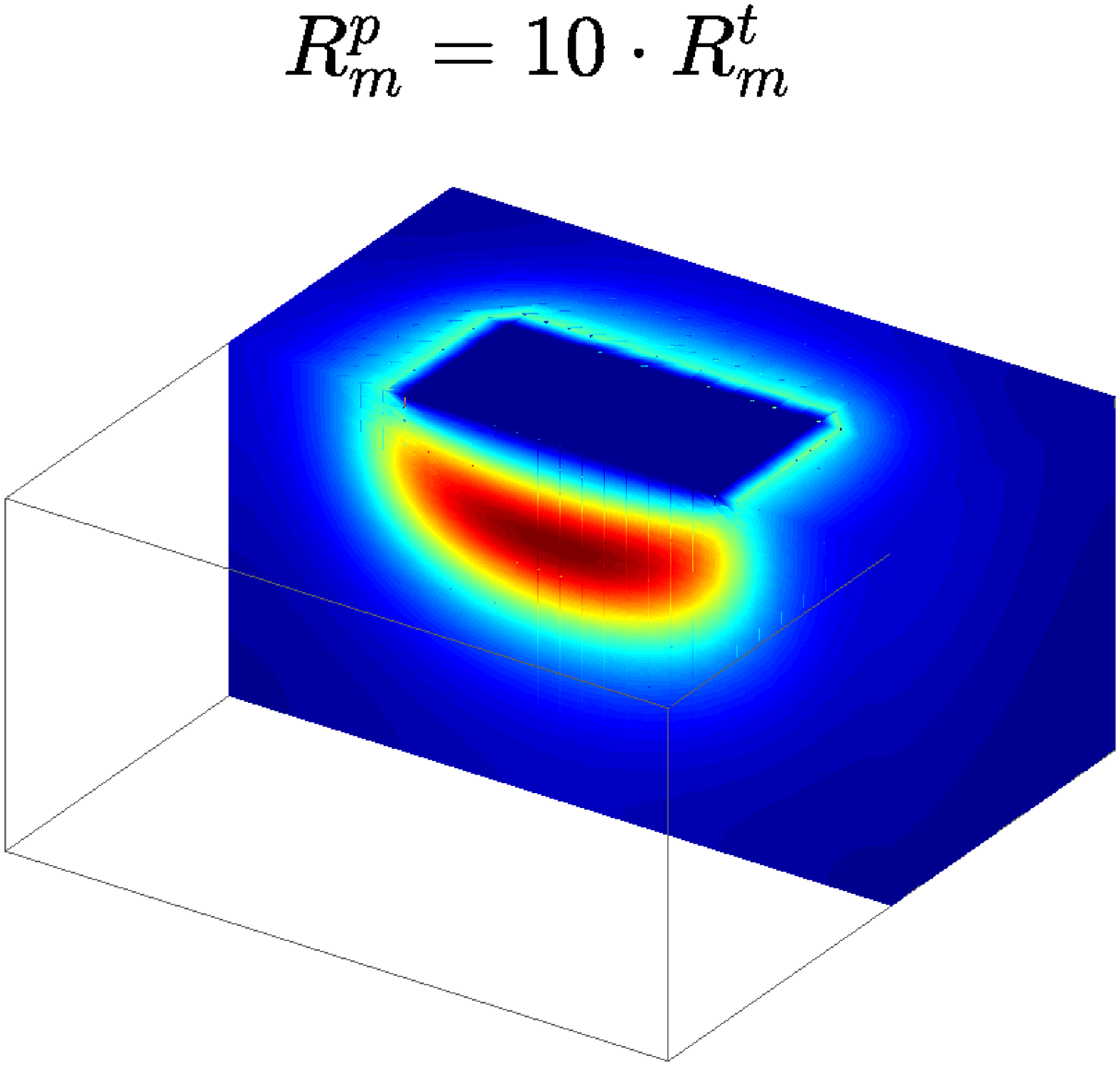}
\caption{Phase-field snapshots for the Mohr-Coulomb crack-driving force for various relations of $R_m^p$ and $R_m^t$.}
\label{fig:BM3D_MC_PF}
\end{figure}

In a next step we want to examine the  Mohr-Coulomb  crack-driving force in more detail. Because this model depends not only on the tensile strength $R_m^t$ but also on the compressive strength $R_m^p$ we investigate different relations between these two parameters. In Fig.~\ref{fig:BM3D_MC}   different load-displacement curves are displayed. The curves with $R_m^t>R_m^p$ have a very similar shape; the graph with $R_m^t<R_m^p$ differs. This is confirmed by the phase-field snapshots of Fig.~\ref{fig:BM3D_MC_PF} which shows that only for  $R_m^t>R_m^p$ the crack nucleates inside the bloc and leads to the specific shell-like breakage.  For $R_m^t<R_m^p$ the crack starts growing at the upper surface. These different mechanisms are caused by the reverse influence of the compressive and tensile states in the driving force.
%

\subsubsection{Finite deformation analysis}
\begin{figure}
\centering
\psfrag{D}[cc]{\small{$\bar u$ [mm]}}
\psfrag{F}[cc][cr]{\small{$F$ [N]}}
\includegraphics[width=0.7\textwidth]{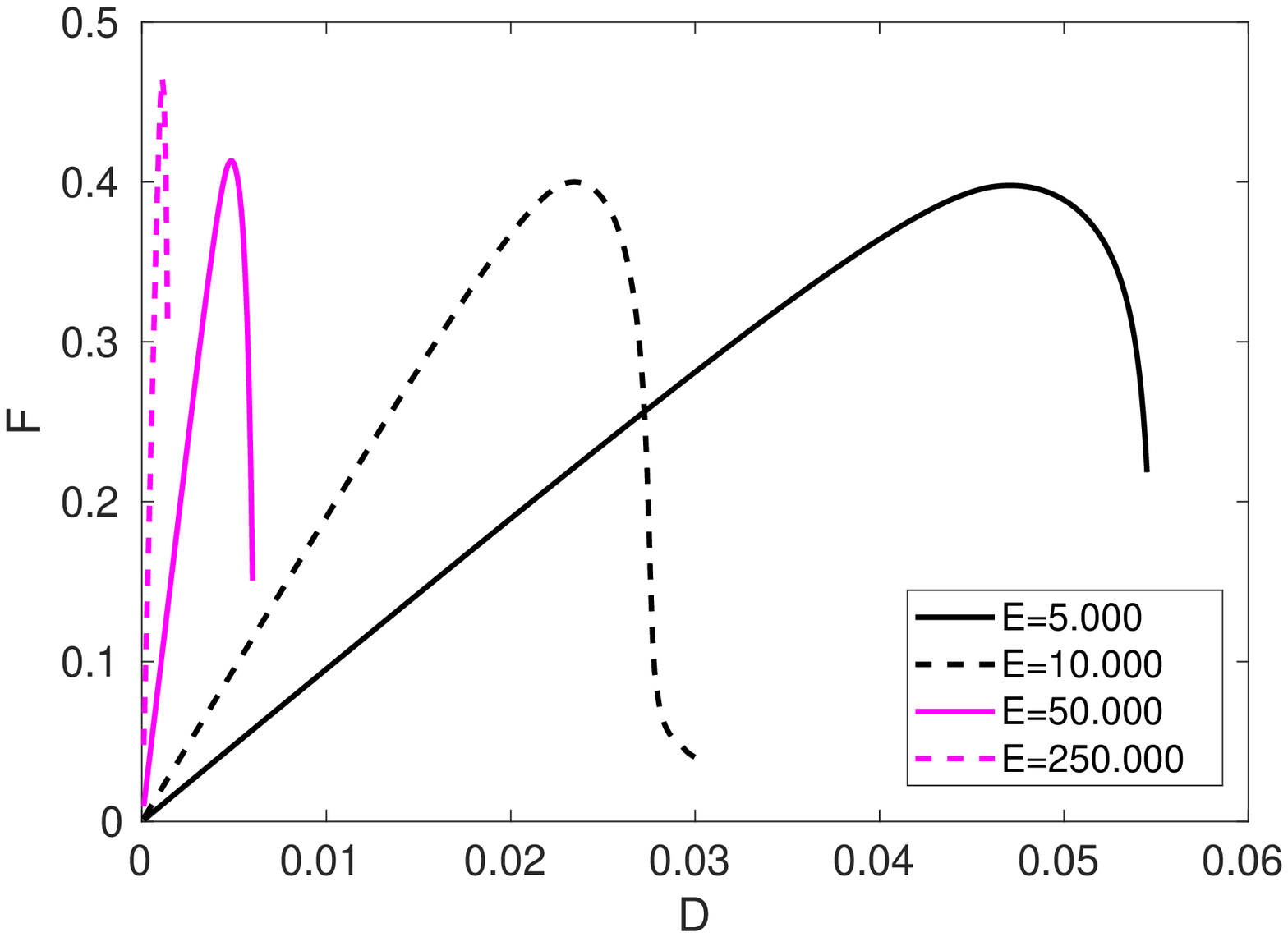}
\caption{Load-displacement curves of the conchoidal fracture example with different elastic modulus $E$.}
\label{fig:BM3D_varyE}
\end{figure}

\begin{figure} 
\centering
\includegraphics[width=0.6\textwidth]{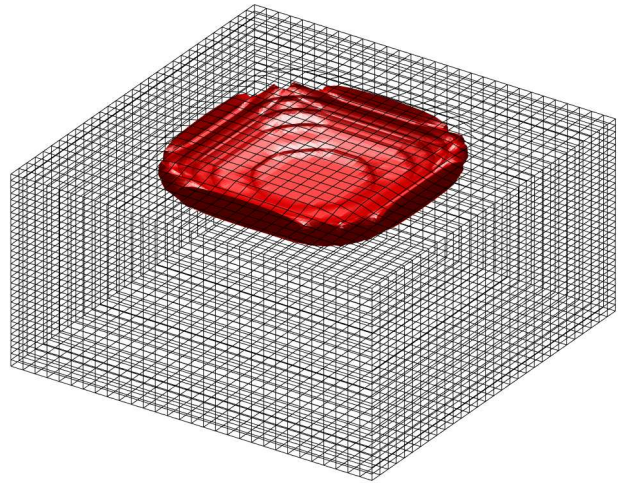}
\caption{Conchoidal fracture surface as a result of rupture in the block with $E=10 000\,\mathrm{N/mm^2}$.}
\label{fig:NHisosurface}
\end{figure}

Finally, the proposed theory is applied in finite elasticity by decreasing the elastic modulus $E$ of this conchoidal fracture example down to a rubbery behavior of the block. The remaining  properties of the material are kept constant as $\nu = 0.25$ and the Griffith energy as $\mathcal{G}_c = 1$ N/mm and follow now a Neo-Hookean material model with energy density (\ref{PsiAdditiveSplitDevVol}-\ref{MooneyRivlinNum}) and $k=0$. The results of the phase-field fracture computations are displayed in Fig.~\ref{fig:BM3D_varyE} and \ref{fig:NHisosurface}. All computations result in the typical conchoidal crack surface which is initiated inside the bloc. The   load-displacement curves show that with lower elastic modulus $E$ the prescribed displacement at crack initiation is higher. The force, however, is almost constant and determined by the material parameter.

We remark, that  everybody may observe  such conch-like fracture surfaces of soft material  in 'ruptured' Swiss cheese for example.


\section{Concluding remarks}\label{sec:conclusion}

The popular phase-field model of fracture in quasi-brittle materials is based on a variational  approach of energy optimization. It has a sound mathematical basis
and   is known to converge for regularization length $l_c\to 0$ to   Griffith' sharp interface model of brittle fracture. For general
numerical computations, however, the crack-driving force cannot follow from the same potential then the deformation field. The physical
restriction that crack growth  presumes a state of tension, requires a split of the body's energy into tensional, crack-driving components
and  insensitive compressive   remainders. In consequence, the energy minimization needs to be performed in a hybrid variational form  with different potentials for phase-field and deformation.

This hybrid form is neither unique 
nor obvious from the physical point of view. Therefore we suggest here, to drive the evolution of the phase-field
without energy minimization, using  physically motivated and ad-hoc formulated driving forces instead. 
We introduce different versions of crack-driving forces, 
compare their pros and cons and perform several numerical computations to illustrate their applicability. In particular we point out
that practical examples, like the Brazilian   test of concrete and the conchoidal fracture of glassy material,
require a crack-driving force specifically derived from fracture mechanics.

%

\section*{Acknowledgements}
The authors  gratefully acknowledge the support of the Deutsche Forschungsgemeinschaft (DFG) under the project "Large-scale simulation of pneumatic and hydraulic fracture with a phase-field approach"  as part of the Priority Program SPP 1748 "Reliable simulation techniques in solid mechanics. Development of non-standard discretization methods, mechanical and mathematical analysis".

\section*{Appendix}\label{sec:appendix}
To find a relation between the Griffith energy $\mathcal{G}_c$ and the fracture strain $\epsilon_c$ 
in phase-field fracture we consider a crack in one dimension. With  energy $\Psi^e_0=\frac12 E \epsilon ^2$, degradation function $g(z)=(1-z)^2$ and crack-surface density function \reff{gamma2O} the variation of energy 
gives
\begin{align*}
    2(1-z) \Psi^e_0 - \frac{\mathcal{G}_c}{l_c}\left( z - l_c^2 \triangle z \right) = 0 \,.
\end{align*}
Now we neglect the regularizing terms to obtain
\begin{align*}
    \left(2 \frac{l_c}{\mathcal{G}_c}\Psi^e_0 +1\right)(1-z) = 1
\end{align*}
Resolving for $1-z$ 
and with the degradated form of \reff{stressesdWdeps} 
it follows a non-monotonous function for the degrading stresses
\begin{align*}
    \sigma = g(z) E \epsilon  =  \left( \frac{l_c}{\mathcal{G}_c}E \epsilon^2 +1\right)^{-2}E \epsilon  \,.
\end{align*}
Its maximum value defines the critical fracture strain $\epsilon_c$.
\begin{align*}
    \epsilon_c= \sqrt{\frac{\mathcal{G}_c}{3 l_c E}}
\end{align*}


\begin{thebibliography}{45}
\providecommand{\natexlab}[1]{#1}
\providecommand{\url}[1]{\texttt{#1}}
\expandafter\ifx\csname urlstyle\endcsname\relax
  \providecommand{\doi}[1]{doi: #1}\else
  \providecommand{\doi}{doi: \begingroup \urlstyle{rm}\Url}\fi

\bibitem[Fro()]{Frocht1948}


\bibitem[Ambati et~al.(2015)Ambati, Gerasimov, and
  De~Lorenzis]{ambati2015phase}
M.~Ambati, T.~Gerasimov, and L.~De~Lorenzis.
\newblock Phase-field modeling of ductile fracture.
\newblock \emph{Computational Mechanics}, 55\penalty0 (5):\penalty0 1017--1040,
  2015.

\bibitem[Ambati et~al.(2016)Ambati, Kruse, and De~Lorenzis]{ambati2016phase}
M.~Ambati, R.~Kruse, and L.~De~Lorenzis.
\newblock A phase-field model for ductile fracture at finite strains and its
  experimental verification.
\newblock \emph{Computational Mechanics}, 57\penalty0 (1):\penalty0 149--167,
  2016.

\bibitem[Ambrosio and Tortorelli(1990)]{AmbrosiTortorelli}
L.~Ambrosio and V.~M. Tortorelli.
\newblock Approximation of functionals depending on jumps by elliptic
  functionals via {$\Gamma$}-convergence.
\newblock \emph{Communications on Pure and Applied Mathematics}, 43\penalty0
  (8):\penalty0 999--1036, 1990.
\newblock ISSN 1097-0312.
\newblock \doi{10.1002/cpa.3160430805}.
\newblock URL \url{http://dx.doi.org/10.1002/cpa.3160430805}.

\bibitem[Amestoy and Leblond(1992)]{Amestoy1992}
M.~Amestoy and J.~Leblond.
\newblock Crack paths in plane situations - ii. detailed form of the expansion
  of the stress intensity factors.
\newblock \emph{International Journal of Solids and Structures}, 29:\penalty0
  465--501, 1992.

\bibitem[Amor et~al.(2009)Amor, Marigo, and Maurini]{amor2009regularized}
H.~Amor, J.-J. Marigo, and C.~Maurini.
\newblock Regularized formulation of the variational brittle fracture with
  unilateral contact: Numerical experiments.
\newblock \emph{Journal of the Mechanics and Physics of Solids}, 57\penalty0
  (8):\penalty0 1209--1229, 2009.

\bibitem[Bilby and Cardew(1975)]{Bilby1975}
B.~Bilby and G.~Cardew.
\newblock The crack with a kinked tip.
\newblock \emph{International Journal of Fracture}, 11\penalty0 (4):\penalty0
  708--712, 1975.

\bibitem[Bilgen et~al.(2017)Bilgen, Kopani{\v c}{\'a}kov{\'a}, Krause, and
  Weinberg]{Bilgen_etal2017}
C.~Bilgen, A.~Kopani{\v c}{\'a}kov{\'a}, R.~Krause, and K.~Weinberg.
\newblock A phase-field approach to conchoidal fracture.
\newblock \emph{Meccanica}, pages 1--17, 2017.

\bibitem[Bleyer and Alessi(2018)]{bleyer2018phase}
J.~Bleyer and R.~Alessi.
\newblock Phase-field modeling of anisotropic brittle fracture including
  several damage mechanisms.
\newblock \emph{Computer Methods in Applied Mechanics and Engineering},
  336:\penalty0 213--236, 2018.

\bibitem[Borden et~al.(2012)Borden, Verhoosel, Scott, Hughes, and
  Landis]{Borden_etal2012}
M.~Borden, C.~Verhoosel, M.~Scott, T.~Hughes, and C.~Landis.
\newblock A phase-field description of dynamic brittle fracture.
\newblock \emph{Comput. Methods Appl. Mech. Engrg.}, 217--220:\penalty0 77--95,
  2012.

\bibitem[Borden et~al.(2014)Borden, Hughes, Landis, and
  Verhoosel]{Borden_etal2014}
M.~J. Borden, T.~J. Hughes, C.~M. Landis, and C.~V. Verhoosel.
\newblock A higher-order phase-field model for brittle fracture: Formulation
  and analysis within the isogeometric analysis framework.
\newblock \emph{Computer Methods in Applied Mechanics and Engineering},
  273\penalty0 (0):\penalty0 100 -- 118, 2014.
\newblock ISSN 0045-7825.
\newblock \doi{http://dx.doi.org/10.1016/j.cma.2014.01.016}.
\newblock URL
  \url{http://www.sciencedirect.com/science/article/pii/S0045782514000292}.

\bibitem[Bourdin et~al.(2008)Bourdin, Francfort, and
  Marigo]{BourdinFrancfortMarigo2008}
B.~Bourdin, G.~Francfort, and J.-J. Marigo.
\newblock \emph{The Variational Approach to Fracture}.
\newblock Springer-Verlag, 2008.

\bibitem[Chambolle et~al.(2009)Chambolle, Francfort, and Marigo]{Chambolle2009}
A.~Chambolle, G.~A. Francfort, and J.-J. Marigo.
\newblock When and how do cracks propagate?
\newblock \emph{Journal of the Mechanics and Physics of Solids}, 57\penalty0
  (9):\penalty0 1614--1622, 2009.

\bibitem[Dally and Weinberg(2017)]{WeinbergDally2015}
T.~Dally and K.~Weinberg.
\newblock The phase-field approach as a tool for experimental validations in
  fracture mechanics.
\newblock \emph{Continuum Mechanics and Thermodynamics}, 29\penalty0
  (4):\penalty0 947--956, 2017.

\bibitem[Doll and Schweizerhof(2000)]{DollSchweizerhof2000}
S.~Doll and K.~Schweizerhof.
\newblock On the development of volumetric strain energy functions.
\newblock \emph{Journal of Applied Mechanics}, 67:\penalty0 17--20, 2000.

\bibitem[Ehlers and Luo(2017)]{ehlers2017phase}
W.~Ehlers and C.~Luo.
\newblock A phase-field approach embedded in the theory of porous media for the
  description of dynamic hydraulic fracturing.
\newblock \emph{Computer Methods in Applied Mechanics and Engineering},
  315:\penalty0 348--368, 2017.

\bibitem[Erdogan and Sih(1963)]{Erdogan1963}
F.~Erdogan and G.~Sih.
\newblock On the crack extension in plates under plane loading and transverse
  shear.
\newblock \emph{Journal of basic engineering}, 85\penalty0 (4):\penalty0
  519--525, 1963.

\bibitem[Francfort and Marigo(1998)]{FrancfortMarigo1998}
G.~Francfort and J.-J. Marigo.
\newblock Revisiting brittle fracture as an energy minimization problem.
\newblock \emph{Journal of the Mechanics and Physics of Solids}, 46:\penalty0
  1319--1342, 1998.
\newblock \doi{10.1016/S0022-5096(98)00034-9}.
\newblock URL \url{http://dx.doi.org/10.1016/S0022-5096(98)00034-9}.

\bibitem[Freddi and Royer-Carfagni(2010)]{FreddyCarfagni2010}
F.~Freddi and G.~Royer-Carfagni.
\newblock Regularized variational theories of fracture: a unified approach.
\newblock \emph{Journal of the Mechanics and Physics of Solids}, 58\penalty0
  (8):\penalty0 1154--1174, 2010.

\bibitem[Goldstein and Salganik(1974)]{Goldstein1974}
R.~Goldstein and R.~Salganik.
\newblock Brittle fracture of solids with arbitrary cracks.
\newblock \emph{\emph{{International Journal of Fracture}}}, 10:\penalty0
  507--523, 1974.

\bibitem[Gross and Selig(2011)]{GrossSelig2011}
D.~Gross and T.~Selig.
\newblock \emph{Bruchmechanik - Mit einer Einführung in die Mikromechanik.}
\newblock Springer, 5. auflage edition, 2011.

\bibitem[Hartmann and Neff(2003)]{HarNef03PGPH}
S.~Hartmann and P.~Neff.
\newblock Polyconvexity of generalized polynomial-type hyperelastic strain
  energy functions for near-incompressibility.
\newblock \emph{Int.\ J.\ Solids and Structures}, 40:\penalty0 2767--2791,
  2003.

\bibitem[Heider and Markert(2017)]{heider2017phase}
Y.~Heider and B.~Markert.
\newblock A phase-field modeling approach of hydraulic fracture in saturated
  porous media.
\newblock \emph{Mechanics Research Communications}, 80:\penalty0 38--46, 2017.

\bibitem[Henry and Levine(2004)]{HenryLevine2004}
H.~Henry and H.~Levine.
\newblock Dynamic instabilities of fracture under biaxial strain using a phase
  field model.
\newblock \emph{Physics Review Letters}, 93:\penalty0 105505, 2004.

\bibitem[Hesch et~al.(2016)Hesch, Schu{\ss}, Dittmann, Franke, and
  Weinberg]{Hesch_etmany2015JCP}
C.~Hesch, S.~Schu{\ss}, M.~Dittmann, M.~Franke, and K.~Weinberg.
\newblock Isogeometric analysis and hierarchical refinement for higher-order
  phase-field models.
\newblock \emph{Computer Methods in Applied Mechanics and Engineering},
  303:\penalty0 185--207, 2016.

\bibitem[Hesch et~al.(2017)Hesch, Gil, Ortigosa, Dittmann, Bilgen, Betsch,
  Franke, Janz, and Weinberg]{Hesch_etmany2017}
C.~Hesch, A.~J. Gil, R.~Ortigosa, M.~Dittmann, C.~Bilgen, P.~Betsch, M.~Franke,
  A.~Janz, and K.~Weinberg.
\newblock A framework for polyconvex large strain phase-field methods to
  fracture.
\newblock \emph{Computer Methods in Applied Mechanics and Engineering},
  317:\penalty0 649--683, 2017.

\bibitem[Hodgdon and Sethna(1993)]{Hodgdon1993}
J.~Hodgdon and J.~Sethna.
\newblock Derivation of a general threedimensional crack-propagation law - a
  generalization of the principle of local symmetry.
\newblock \emph{\emph{{Physical Review}}}, 47:\penalty0 4831--4840, 1993.

\bibitem[Karma et~al.(2001)Karma, Kessler, and Levine]{Karma_etal2001}
A.~Karma, D.~Kessler, and H.~Levine.
\newblock Phase-field model of mode iii dynamic fracture.
\newblock \emph{Physical Review Letter}, 81:\penalty0 045501, 2001.

\bibitem[Katzav et~al.(2007)Katzav, Adda-Bedia, and
  Derrida]{katzav2007fracture}
E.~Katzav, M.~Adda-Bedia, and B.~Derrida.
\newblock Fracture surfaces of heterogeneous materials: a 2d solvable model.
\newblock \emph{EPL (Europhysics Letters)}, 78\penalty0 (4):\penalty0 46006,
  2007.

\bibitem[Khosravani et~al.(2018)Khosravani, Silani, and
  Weinberg]{Khosravanietal2018TAFM}
M.~R. Khosravani, M.~Silani, and K.~Weinberg.
\newblock Fracture studies of ultra-high performance concrete using dynamic
  brazilian tests.
\newblock \emph{Theoretical and Applied Fracture Mechanics}, 93:\penalty0
  302--310, 2018.

\bibitem[Kuhn et~al.(2016)Kuhn, Noll, and M{\"u}ller]{kuhn2016phase}
C.~Kuhn, T.~Noll, and R.~M{\"u}ller.
\newblock On phase field modeling of ductile fracture.
\newblock \emph{GAMM-Mitteilungen}, 39\penalty0 (1):\penalty0 35--54, 2016.

\bibitem[Leblond(1989)]{Leblond1989}
J.~Leblond.
\newblock Cracks paths in plane situations - i. general form of the expansion
  of the stress intensity factors.
\newblock \emph{\emph{{International Journal of Solids and Structures}}},
  362:\penalty0 295--296, 1989.

\bibitem[Li et~al.(2016)Li, Marigo, Guilbaud, and Potapov]{LiMarigo2017}
T.~Li, J.-J. Marigo, D.~Guilbaud, and S.~Potapov.
\newblock Gradient damage modeling of brittle fracture in an explicit dynamics
  context.
\newblock \emph{International Journal for Numerical Methods in Engineering},
  108\penalty0 (11):\penalty0 1381--1405, 2016.

\bibitem[Miehe and Sch{\"a}nzel(2014)]{miehe2014phase}
C.~Miehe and L.-M. Sch{\"a}nzel.
\newblock Phase field modeling of fracture in rubbery polymers. part i: Finite
  elasticity coupled with brittle failure.
\newblock \emph{Journal of the Mechanics and Physics of Solids}, 65:\penalty0
  93--113, 2014.

\bibitem[Miehe et~al.(2010)Miehe, Hofacker, and Welschinger]{miehe2010phase}
C.~Miehe, M.~Hofacker, and F.~Welschinger.
\newblock A phase field model for rate-independent crack propagation: {R}obust
  algorithmic implementation based on operator splits.
\newblock \emph{Comput. Methods Appl. Mech. Engrg.}, 199:\penalty0 2765--2778,
  2010.

\bibitem[Miehe et~al.(2015{\natexlab{a}})Miehe, Hofacker, Sch{\"a}nzel, and
  Aldakheel]{miehe2015phase}
C.~Miehe, M.~Hofacker, L.-M. Sch{\"a}nzel, and F.~Aldakheel.
\newblock Phase field modeling of fracture in multi-physics problems. part ii.
  coupled brittle-to-ductile failure criteria and crack propagation in
  thermo-elastic-plastic solids.
\newblock \emph{Computer Methods in Applied Mechanics and Engineering},
  294:\penalty0 486--522, 2015{\natexlab{a}}.

\bibitem[Miehe et~al.(2015{\natexlab{b}})Miehe, Schänzel, and
  Ulmer]{miehe_etal_2015a}
C.~Miehe, L.-M. Schänzel, and H.~Ulmer.
\newblock Phase-field modeling of fracture in multi-physics problems. part i.
  balance of crack surface and failure criteria for brittle crack propagation
  in thermo-elasitc solids.
\newblock \emph{Comput. Methods Appl. Mech. Engrg.}, 294:\penalty0 449--485,
  2015{\natexlab{b}}.

\bibitem[Negri(2013)]{negri2013phase}
M.~Negri.
\newblock From phase field to sharp cracks: Convergence of quasi-static
  evolutions in a special setting.
\newblock \emph{Applied Mathematics Letters}, 26\penalty0 (2):\penalty0
  219--224, 2013.

\bibitem[Sargado et~al.(2018)Sargado, Keilegavlen, Berre, and
  Nordbotten]{sargado2018high}
J.~M. Sargado, E.~Keilegavlen, I.~Berre, and J.~M. Nordbotten.
\newblock High-accuracy phase-field models for brittle fracture based on a new
  family of degradation functions.
\newblock \emph{Journal of the Mechanics and Physics of Solids}, 111:\penalty0
  458--489, 2018.

\bibitem[Teichtmeister et~al.(2017)Teichtmeister, Kienle, Aldakheel, and
  Keip]{teichtmeister2017phase}
S.~Teichtmeister, D.~Kienle, F.~Aldakheel, and M.-A. Keip.
\newblock Phase field modeling of fracture in anisotropic brittle solids.
\newblock \emph{International Journal of Non-Linear Mechanics}, 97:\penalty0
  1--21, 2017.

\bibitem[Thomas et~al.(2018)Thomas, Bilgen, and
  Weinberg]{Thomas_etal_2018_GAMM}
M.~Thomas, C.~Bilgen, and K.~Weinberg.
\newblock Phase‐field fracture at finite strains based on modified
  invariants: A note on its analysis and simulations.
\newblock \emph{GAMM-Mitteilungen}, 40 (3):\penalty0 207--237, 2018.
\newblock \doi{10.1002/gamm.201730004}.

\bibitem[Verhoosel and de~Borst(2013)]{deborst2013}
C.~Verhoosel and R.~de~Borst.
\newblock A phase-field model for cohesive fracture.
\newblock \emph{Int. J. Numer. Methods Eng.}, 2013.

\bibitem[Wallner(1939)]{Wallner1939}
H.~Wallner.
\newblock Linienstrukturen an bruchfl{\"a}chen.
\newblock \emph{Zeitschrift f{\"u}r Physik}, 114\penalty0 (5-6):\penalty0
  368--378, 1939.

\bibitem[Weinberg and Hesch(2015)]{WeinbergHesch2015}
K.~Weinberg and C.~Hesch.
\newblock A high-order finite-deformation phase-field approach to fracture.
\newblock \emph{Continuum Mech. Thermodyn.}, 2015.
\newblock \doi{10.1007/s00161-015-0440-7}.

\bibitem[Weinberg et~al.(2016)Weinberg, Dally, Schu{\ss}, Werner, and
  Bilgen]{WeinbergGAMMMitteilungen2016}
K.~Weinberg, T.~Dally, S.~Schu{\ss}, M.~Werner, and C.~Bilgen.
\newblock Modeling and numerical simulation of crack growth and damage with a
  phase field approach.
\newblock \emph{GAMM-Mitteilungen}, 39(1):\penalty0 55--77, 2016.

\end{thebibliography}

\end{document}